\documentclass[parskip=half,bibliography=totoc]{scrartcl}

\usepackage[automark]{scrlayer-scrpage}								
\usepackage{amsfonts}												
\usepackage{amsmath}												
\usepackage{mathtools}												
\usepackage{stmaryrd}												
\usepackage{amssymb}												
\usepackage{mathrsfs}												
\usepackage{accents}												
\usepackage[T1]{fontenc}											
\usepackage[utf8]{inputenc}											
\usepackage[top=1in,bottom=1in,left=1.25in,right=1.25in]{geometry}	
\usepackage{enumerate} 												
\usepackage{authblk}												
\usepackage{abstract}												
\usepackage{etoolbox}												

\usepackage{tikz-cd}												

\usepackage{amsthm}													

\usepackage[colorlinks=true,citecolor=blue]{hyperref}				
\usepackage[noabbrev]{cleveref}										

\newcommand{\pd}[2]{\frac{\partial #1}{\partial #2}}
\renewcommand{\d}{\mathrm{d}}

\newcommand{\HP}{\mathbb{H}\mathrm{P}}

\newcommand{\R}{\mathbb{R}}
\newcommand{\C}{\mathbb{C}}
\renewcommand{\H}{\mathbb{H}}

\DeclareMathOperator{\tr}{tr}

\DeclareMathOperator{\End}{End}

\DeclareMathOperator{\id}{id}

\DeclareMathOperator{\scal}{scal}

\newcommand{\abs}[1]{\lvert #1 \rvert}
\newcommand{\norm}[1]{\lVert #1 \rVert}

\newcommand{\mf}[1]{\mathfrak{#1}}
\newcommand{\mc}[1]{\mathcal{#1}}

\newcommand{\h}{\mathrm{H}}
\newcommand{\q}{\mathrm{Q}}

\newtheoremstyle{mythm}
{}
{}
{\slshape}
{}
{\bfseries\sffamily}
{.}
{ }
{}
\newtheoremstyle{mydef}
{}
{}
{}
{}
{\bfseries\sffamily}
{.}
{ }
{}

\theoremstyle{mythm}
\newtheorem{thm}{Theorem}[section]

\newtheorem{lem}[thm]{Lemma}
\theoremstyle{mydef}
\newtheorem{mydef}[thm]{Definition}
\newtheorem{rem}[thm]{Remark}

\newenvironment{myproof}[1][\proofname]{
	\proof[\sffamily\upshape#1]
}{\endproof}

\newcommand{\proofclear}{\hfill \qedsymbol}

\clearscrheadfoot
\ihead[]{\headmark}
\ohead[]{\pagemark}
\deffootnote[1em]{0em}{1em}{%
	\textsuperscript{\thefootnotemark}%
}
\setfootnoterule{3em}

\newcommand\numberthis{\stepcounter{equation}\tag{\theequation}}

\newenvironment{numberedlist}{\begin{enumerate}[\upshape(i)]}{\end{enumerate}}

\apptocmd{\sloppy}{\hbadness 10000\relax}{}{}

\title{Curvature of quaternionic Kähler manifolds with $S^1$-symmetry}
\author{V.\ Cortés}
\author{A.\ Saha}
\author{D.\ Thung}
\affil{\normalsize  Department of Mathematics \\
University of Hamburg\\
Bundesstra\ss e 55, D-20146 Hamburg, Germany}

\date{}

\begin{document}
\maketitle

\begin{abstract}
	We study the behavior of connections and curvature under the HK/QK correspondence, proving simple formulae expressing the Levi-Civita connection and Riemann curvature tensor on the quaternionic K\"ahler side in terms of the initial hyper-K\"ahler data. Our curvature formula refines a well-known decomposition theorem due to Alekseevsky. As an application, we compute the norm of the curvature tensor for a series of complete quaternionic K\"ahler manifolds arising from flat hyper-K\"ahler manifolds. We use this to deduce that these manifolds are of cohomogeneity one.\par
	\emph{Keywords: quaternionic K\"ahler manifolds, HK/QK correspondence, $c$-map, one-loop deformation, curvature, cohomogeneity one}\par
	\emph{MSC classification: 53C26.}
\end{abstract}

\clearpage

\tableofcontents
\clearpage

\section{Introduction}

This paper concerns itself with the hyper-K\"ahler/quaternionic K\"ahler (HK/QK) correspondence, which is a duality between hyper-K\"ahler and quaternionic K\"ahler manifolds with (infinitesimal) circle symmetries. It was discovered by Haydys \cite{Hay2008}, and has since attracted the attention of physicists and mathematicians alike \cite{APP2011,ACM2013,Hit2013,ACDM2015,MS2015}. 

One of the most important applications of the HK/QK correspondence is to the study of a construction of quaternionic K\"ahler manifolds known as the (supergravity) $c$-map \cite{FS1990,Hit2009}. This construction, first discovered by theoretical physicists, was extended to produce one-parameter families of quaternionic K\"ahler metrics \cite{RSV2006}, known as one-loop deformed Ferrara--Sabharwal metrics (reflecting their physical origins).

The starting point for the $c$-map construction is not a hyper-K\"ahler manifold, but rather a projective special K\"ahler (PSK) manifold $\bar M$, which is by definition the K\"ahler quotient of a corresponding conical affine special K\"ahler (CASK) manifold $M$ with respect to a preferred circle action. The $c$-map and its one-loop deformation are most conveniently described by first passing to $M$, then to its cotangent bundle $N=T^*M$, which comes with a natural (pseudo-)hyper-K\"ahler structure. The circle action on $M$ moreover lifts to a natural circle symmetry of $N$, making it possible to apply the HK/QK correspondence to $N$ to produce a quaternionic K\"ahler manifold $\bar N$ of the same dimension as $N$ (see~\cite{APP2011,ACM2013,ACDM2015} for details). The resulting quaternionic K\"ahler metric admits a (non-zero) Killing field, and in fact any quaternionic K\"ahler manifold with a Killing field can be locally obtained in this fashion.

An important feature of the HK/QK correspondence is that, in its applications to the $c$-map, the one-parameter family of one-loop deformed Ferrara--Sabharwal metrics $g_\text{FS}^c$ can be studied in uniform fashion. This is due to the fact that, on the hyper-K\"ahler side, the deformation parameter $c$ corresponds to the additive freedom in choosing a Hamiltonian for the circle action, which is easy to handle.

In this paper, we study the behavior of curvature tensors under the HK/QK correspondence. We make use of the point of view put forward by Macia and Swann \cite{MS2015}, which builds on Swann's twist construction \cite{Swa2010}, casting the HK/QK correspondence as a variation on this general method. The main advantage of this approach is that the twist construction gives explicit relations between certain tensor fields on both sides of the HK/QK correspondence. 

In \Cref{sec:background}, we briefly recall the general twist formalism and its application to the HK/QK correspondence. Using these methods, we investigate how connections and curvature tensors change under the HK/QK correspondence in \Cref{sec:maintheorems}. This culminates in elegant expressions for the Levi-Civita connection and curvature of the quaternionic K\"ahler metric in terms of the hyper-K\"ahler data on the other side of the correspondence. Our formulas apply in particular to the quaternionic K\"ahler manifolds that arise from the $c$-map, whose curvature tensors have previously been studied in the physics literature \cite{dWVvP1993} (see also \cite{CDJL2017} for curvature formulas in the special case of the so-called $q$-map).

In \Cref{sec:applications}, we apply our results to a series of examples considered in~\cite{CST2020}. There, we showed that, under the (one-loop deformed) $c$-map, automorphisms of the initial PSK manifold yield isometries of the resulting quaternionic K\"ahler metric. This leads to lower bounds on the degree of symmetry of one-loop deformed $c$-map spaces. In particular, applying the one-loop deformed $c$-map to a homogeneous PSK manifold results in a quaternionic K\"ahler metric of cohomogeneity at most one. By contrast, our curvature formula can be used to give upper bounds on the degree of symmetry. To demonstrate this, we consider a one-parameter deformation of the symmetric spaces $SU(n+1,2)/S(U(n+1)\times U(2))$, which arise by applying the $c$-map to complex hyperbolic spaces, regarded as homogeneous PSK manifolds. Our curvature formula enables us to compute the norm of the curvature endomorphisms of these quaternionic K\"ahler metrics, which the HK/QK correspondence relates to a series of flat hyper-K\"ahler manifolds. This leads to a proof that these quaternionic K\"ahler manifolds are not locally homogeneous. In fact, combining this with our results in \cite{CST2020}, we establish that the isometry group acts with cohomogeneity precisely one in this family of examples.

{\bfseries Acknowledgements}

This work was supported by the German Science Foundation (DFG) under the Research Training Group 1670, and under Germany’s Excellence Strategy – EXC 2121 ``Quantum Universe'' – 390833306.

\clearpage

\section{Twist construction and HK/QK correspondence}
\label{sec:background}

In this section, we describe the HK/QK correspondence as a variation on the twist construction \cite{Swa2010}, following Macia and Swann \cite{MS2015}. For a more detailed review, we refer the reader to \cite{CST2020}. 

\subsection{Twists}

The twist construction is a general method to produce new manifolds with (infinitesimal) torus actions out of others. In the following, we will restrict to free circle actions, since this is the case of interest for the HK/QK correspondence. Let $N$ be a manifold endowed with a free circle action, generated by a vector field $Z$. Our goal is to construct another manifold, which carries a ``twisted'' circle action. 

To this end, we first pick a closed, integral two-form $\omega$ on $N$ and construct a circle bundle $\pi_N:P\to N$ equipped with a connection $\eta$ with curvature $\omega$. Now we ask for a lift of $Z$ to $P$ which preserves $\eta$ and commutes with $X_P$, the generator of the principal circle action on $P$. Such a lift exists if and only if $Z$ is $\omega$-Hamiltonian, i.e.~there exists a function $f\in C^\infty(N)$ such that $\iota_Z\omega=-\d f$. In this case, it is given by $Z_P=\tilde Z+ \pi_N^* f X_P$, where $\tilde Z$ denotes the $\eta$-horizontal lift of $Z$. We will always assume that $f$ is nowhere-vanishing, so that $Z_P$ is transversal to the horizontal distribution $\mc H\coloneqq \ker \eta\subset TP$. Now we define the twist of $N$ with respect to the twist data $(Z,\omega,f)$ as $\bar N\coloneqq X/\langle Z_P\rangle$. Then $P$ has the structure of a double fibration:
\begin{equation*}
	\begin{tikzcd}
		N & P \ar[l,"\pi_N"'] \ar[r,"\pi_{\bar N}"] & \bar N
	\end{tikzcd}
\end{equation*}
and the action of $X_P$ pushes down to a vector field $Z_{\bar N}$ on $\bar N$ which plays a dual role to $Z$.

Since both $X_P$ and $Z_P$ are transversal to $\mc H$, we have pointwise identifications of the tangent spaces of both $N$ and $\bar N$ with the horizontal subspaces of $TP$. This sets up a correspondence between tensor fields. On each space, we have a notion of horizontal lift to $P$ for arbitrary tensor fields, induced by the usual horizontal lift on vector fields and the composition of pullback and restriction to $\mc H$ for one-forms. We then say that tensor fields $T$ on $N$ and $T'$ on $\bar N$ are $\mc H$-related, denoted $T\sim_{\mc H}T'$, if their horizontal lifts agree. We may also call $T'$ the twist of $T$---the twist data implicitly understood. Note, however, that the existence of a well-defined twist implicitly presumes $Z$-invariance of $T$ (and, dually, $Z_{\bar N}$-invariance of $T'$), so the correspondence only makes sense for $Z$-invariant tensor fields.

This method of carrying invariant tensor fields from $N$ over to $\bar N$, and vice versa, is compatible with tensor products and contractions and respects algebraic relations. However, it typically does not preserve differential conditions, introducing correction terms dependent on the twist data. For example, if $\alpha \sim_{\mc H}\alpha'$ are differential forms, then $\d \alpha-f^{-1}\omega\wedge \iota_Z\alpha \sim_{\mc H}\d \alpha'$.

\subsection{Elementary deformations and HK/QK correspondence}

Let us now consider the case where $(N,g,\omega_1,\omega_2,\omega_3)$ is a (pseudo-)hyper-K\"ahler manifold and comes equipped with a rotating (circle) symmetry, that is, the generating vector field $Z$ of the circle action satisfies $L_Z g=L_Z \omega_1=0$, while $L_Z\omega_2=\omega_3$ and $L_Z\omega_3=-\omega_2$. We also assume that $Z$ is $\omega_1$-Hamiltonian, i.e.~$\iota_Z\omega_1=-\d f_Z$, and that $\omega_1$ is integral, so that $(Z,\omega_1,f_Z)$ can be used as twist data.

The resulting twist manifold will inherit an almost quaternion-Hermitian structure from $N$. However, the fundamental four-form, which is the twist of the $Z$-invariant parallel four-form $\Omega=\sum_{j=1}^3\omega_j\wedge\omega_j$ on $N$, will fail to be parallel because of the general failure of the twist construction to preserve differential relations. This means that the twist manifold is not quaternionic K\"ahler. 

Macia and Swann investigated how to deform the hyper-K\"ahler structure in order to produce a quaternionic K\"ahler manifold out of a hyper-K\"ahler manifold. To this end, they defined the notion of elementary deformations of such manifolds. An elementary deformation $g_\h$ of the hyper-K\"ahler metric $g$ with respect to a rotating symmetry $Z$ (with $\omega_1$-Hamiltonian $f_Z$) is a new metric of the form
\begin{equation*}
g_\h=a \cdot g+b\cdot \bigg((\iota_Z g)^2+\sum_{j=1}^3 (\iota_Z \omega_j)^2\bigg)
\end{equation*}
for smooth, nowhere-vanishing functions $a,b$ on $N$. Note that the second term is equal to $b\cdot g(Z,Z)g|_{\H Z}$, so we think of $g_\h$ as resulting from a conformal scaling composed with an independent scaling along the quaternionic span of $Z$. In \cite{MS2015}, Macia and Swann study the twists of elementary deformations and prove the following result:

\begin{thm}[\cite{MS2015}]\label{thm:HKQK}
	Let $(N,g,\omega_j)$, $j=1,2,3$, be a (pseudo-)hyper-K\"ahler manifold endowed with a rotating circle symmetry with $\iota_Z\omega_1=-\d f_Z$. Then an elementary deformation $g_\h$ of $g$ with respect to $Z$ twists to a (pseudo-)quaternionic K\"ahler metric on the twist manifold $\bar N$ if and only if the twist data are given by $(Z,\omega_\h\coloneqq k(\omega_1+\d\iota_Z g),f_\h\coloneqq k(f_Z+g(Z,Z)))$ and the elementary deformation is
	\begin{equation*}
		g_\h=\frac{B}{f_Z}g+\frac{B}{f_Z^2}\bigg((\iota_Z g)^2+\sum_{j=1}^3 (\iota_Z \omega_j)^2\bigg)
	\end{equation*}
	where $k,B\in \R\setminus\{0\}$ are constants.
\end{thm}

In the following, we will only consider the elementary deformation given in the theorem, and will refer to it as \emph{the} elementary deformation of $g$.

\begin{rem}\leavevmode
	\begin{numberedlist}
		\item The constants $k$ and $B$ simply scale the curvature of $P$ and the metric, respectively. In the following, we will set them to $1$, and assume that $\omega_1$ is integral. Nevertheless, there is a one-parameter freedom in the construction, which arises because we may add a constant to $f_Z$. Thus, one actually obtains a one-parameter family of quaternionic K\"ahler metrics.
		\item In the case of interest for the $c$-map and its one-loop deformation, the resulting quaternionic K\"ahler metric is positive-definite \cite{ACDM2015}.
	\end{numberedlist}
\end{rem}

\subsection{Notation and some useful identities}
\label{sec:notation}

In this section, we introduce some notation and derive a number of identities that will be useful to us in the following computations.

In the following, it will often be convenient to discuss the metric and its K\"ahler forms in uniform manner. Therefore, we introduce $I_0\coloneqq \id$ and $\omega_0\coloneqq g$, so that $\omega_\mu(\cdot,\cdot)=g(I_\mu\cdot,\cdot)$, $\mu=0,1,2,3$. We also set $\alpha_\mu\coloneqq \iota_Z\omega_\mu$. In what follows, Greek indices will be understood to run from $0$ to $3$, while Latin indices run from $1$ to $3$.

We briefly discuss the relation between $g_\h$ and $g$. 

\begin{mydef}\label{def:metriccomparison}
	Define the endomorphism field $\mc K:TN\to TN$ by $g_\h(\mc K X,Y)=g(X,Y)$. Equivalently, regarding both $g_\h$ and $g$ as maps $TN\to T^*N$, we have $\mc K=g_\h^{-1}\circ g$.
\end{mydef}

Since $\sum_\mu \alpha_\mu^2=g(Z,Z)g|_{\H Z}$, we can rewrite the elementary deformation as follows:
\begin{equation*}
	g_\h=\frac{1}{f_Z}g|_{(\H Z)^\perp}+\frac{f_\h}{f_Z^2}g|_{\H Z}
\end{equation*}
This shows that $\mc K|_{(\H Z)^\perp}=f_Z$ and $\mc K|_{\H Z}=\frac{f_Z^2}{f_\h}$.

\begin{lem}\label{lem:metriccomparecommutes}
	$\mc K$ is self-adjoint with respect to $g$ and commutes with each $I_\mu$.
\end{lem}
\begin{myproof}
	The above description of $\mc K$ implies that 
	\begin{equation}\label{eq:metriccompareformula}
		\mc K(X)=f_Z X - \frac{f_Z}{f_\h}\sum_{\lambda=0}^3  \alpha_\lambda(X)I_\lambda Z
	\end{equation}
	Indeed, the right-hand side restricts to the appropriate multiples of the identity on $\H Z$ and its orthogonal complement. Self-adjointness is easy to check. Now we compute
	\begin{align*}
		\mc K(I_\mu X)&=f_Z I_\mu X-\frac{f_Z}{f_\h}\sum_\lambda g(I_\lambda Z,I_\mu X)I_\lambda Z
		=f_Z I_\mu X-\frac{f_Z}{f_\h}\sum_\lambda g(I_\mu^{-1}I_\lambda Z,X)I_\lambda Z\\
		&=f_Z I_\mu X-\frac{f_Z}{f_\h}\sum_\lambda g(I_\lambda Z,X)I_\mu I_\lambda Z=I_\mu \mc K(X)
	\end{align*}
	In passing to the second line, we made the substitution $I_\lambda\mapsto I_\mu I_\lambda$ in the sum. We can do this without introducing any additional signs, because $I_\lambda$ appears twice in the expression.
\end{myproof}

Now we study how $g$ relates to the twisting form $\omega_\h$. 

\begin{lem}
	Denoting the Levi-Civita connection of $g$ by $D$, the one-forms $\alpha_\mu\coloneqq \iota_Z \omega_\mu$ satisfy
	\begin{equation}\label{eq:dalpha}
		\begin{aligned}
			(\d\alpha_0)(A,B)&
			=2g(D_A Z,B)\\
			(\d\alpha_k)(A,B)&=(L_Z\omega_k)(A,B)\\
			&=\omega_k(D_AZ,B)+\omega_k(A,D_B Z)
		\end{aligned}
	\end{equation}
\end{lem}
\begin{myproof}
	The first identity follows from Cartan's formula and the Killing equation, while the second follows by applying Cartan's formula and the fact that $L_Z=D_Z-DZ$ on differential forms, remembering that each $\omega_k$ is parallel.
\end{myproof}

\begin{mydef}
	We define the skew-symmetric (with respect to $g$) endomorphism field $I_\h\coloneqq I_1+2DZ$.
\end{mydef}

By definition of $\omega_\h$ and the first part of the above Lemma, we have $\omega_\h(\cdot,\cdot)=g(I_\h \cdot,\cdot)$.

Using $I_\h$ and the action of $Z$ on $\omega_k$, we deduce the following consequence of \eqref{eq:dalpha}:
\begin{equation*}
	\omega_\mu(D_AZ,B)-\omega_\mu(D_BZ,A)=-\frac{1}{2}(\omega_\mu(I_1 A,B)-\omega_\mu(I_1 B,A))+\delta_{\mu 0}\omega_\h(A,B)
\end{equation*}
where $\delta$ is the Kronecker delta symbol. Now we may sum over $\mu$ to obtain the following identity: 
\begin{equation}\label{eq:sumidentity}
	\begin{aligned}
		&\sum_\mu \big(\omega_\mu(D_AZ,B)-\omega_\mu(D_BZ,A)\big)I_\mu I_1 C\\
		&\quad =-\frac{1}{2}\sum_\mu \big(\omega_\mu(A,B)-\omega_\mu(B,A)\big)I_\mu C + \omega_\h(A,B)I_1 C
	\end{aligned}
\end{equation}
To obtain the right-hand side, we substituted $I_\mu$ for $I_\mu I_1$ in the first terms (the same trick as in the proof of \Cref{lem:metriccomparecommutes}).

\begin{lem}
	$I_\h$ commutes with each $I_\mu$.
\end{lem}
\begin{myproof}
	From $L_Z\omega_1=0$, $L_Z\omega_2=\omega_3$ and $L_Z\omega_3=-\omega_2$ and the second line of \eqref{eq:dalpha}, one obtains
	\begin{align*}
		&2\omega_1(D_A Z,B)+2\omega_1(A,D_BZ)=0=-\omega_1(I_1 A,B)+\omega_1(I_1 B,A)\\\numberthis\label{eq:omegaidentities}
		&2\omega_2(D_A Z,B)+2\omega_2(A,D_BZ)=2\omega_3(A,B)=\omega_1(I_2 A,B)-\omega_1(I_2 B,A)\\
		&2\omega_3(D_A Z,B)+2\omega_3(A,D_BZ)=-2\omega_2(A,B)=\omega_1(I_3 A,B)-\omega_1(I_3 B,A)
	\end{align*}
	Rearranging each line and using our definition of $I_\h$, we can write this compactly as $g(I_k I_\h A,B)=g(I_k I_\h B,A)$ for $k=1,2,3$. We also have $g(I_0 I_\h A,B)=-g(I_0 I_\h B, A)$, and in summary find
	\begin{equation*}
		g(I_\mu I_\h A,B)=-g(I^{-1}_\mu I_\h B,A)
	\end{equation*}
	Now, we conclude with a short computation:
	\begin{align*}
		g(I_\mu I_\h A,B)&=-g(I^{-1}_\mu I_\h B,A)=-g(I_\h B,I_\mu A)=-\omega_\h(B,I_\mu A)=\omega_\h(I_\mu A,B)\\
		&=g(I_\h I_\mu A,B)
	\end{align*}
	which implies the claim.
\end{myproof}

\begin{rem}
	In fact, this line of reasoning shows that, given an arbitrary vector field $Z$ on a hyper-K\"ahler manifold, the condition that $I_1+2DZ$ defines a skew-symmetric endomorphism field which commutes with every $I_\mu$ is equivalent to requiring that $Z$ is an infinitesimal rotating symmetry.
\end{rem}

\subsection{General twisting formulae for connections and curvature}
\label{sec:twistformulae}

Before we can study the behavior of the Riemann curvature tensor under the HK/QK correspondence, we must understand how the Levi-Civita connection changes under elementary deformations and twists.

\begin{lem}
	Let $(N,g,\omega_k)$, $k=1,2,3$ be a hyper-K\"ahler manifold endowed with a rotating circle action, and let $g_\h$ be the elementary deformation of $g$ (cf.~\Cref{thm:HKQK}). If $D$ and $D^\h$ denote the Levi-Civita connections of $g$ and $g_\h$, respectively, we have $D^\h=D+S^\h$, where 
	\begin{equation}\label{eq:S^h}
		2g_\h(S^\h_AB,C)=(D_A g_\h)(B,C)+(D_Bg_\h)(C,A)-(D_Cg_\h)(A,B)
	\end{equation}
	for all vector fields $A,B,C$ on $N$.
\end{lem}
\begin{myproof}
	This follows from cyclically permuting the vector fields $A,B,C$ in the identity
	\begin{equation*}
		0=(D^\h_Ag_\h)(B,C)=(D_Ag_\h)(B,C)-g_\h(S^\h_AB,C)-g_\h(B,S^\h_AC)
	\end{equation*}
	and taking the appropriate signed sum.
\end{myproof}

\begin{mydef}
	Let $N$ be a manifold endowed with twist data, and let $\nabla$ be a $Z$-invariant connection on $N$. We say that $\nabla$ is $\mc H$-related to a connection $\nabla'$ on the twist manifold if, for all $Z$-invariant vector fields $A\sim_{\mc H}A'$, $B\sim_{\mc H}B'$, we have $\nabla'_{A'}B'\sim_{\mc H}\nabla_AB$.
\end{mydef}

\begin{lem}
	Let $(N,g_\h)$ be a Riemannian manifold endowed with the twist data $(Z,\omega_\h,f_\h)$. Assume that $g_\h$ is $Z$-invariant, and let $g_\q$ be the $\mc H$-related metric on the twist manifold $\bar N$. Denote the corresponding Levi-Civita connections by $D^\h$ and $D^\q$. Then $D^\q \sim_{\mc H}D^\h + S^\q$, where
	\begin{equation}\label{eq:S^q}
		2g_\h(S^\q_AB,C)=\frac{1}{f_\h}\Big(g_\h(Z,C)\omega_\h(A,B)-g_\h(Z,A)\omega_\h(B,C)-g_\h(Z,B)\omega_\h(A,C)\Big)
	\end{equation}
	for arbitrary $Z$-invariant vector fields $A,B,C$.
\end{lem}
\begin{myproof}
	We use the Koszul formula for $D^\q$ using vector fields $A',B',C'$, $\mc H$-related to $A,B,C$:
	\begin{align*}
		2g_\q(D^\q_{A'} B',C')&=A'g_\q(B',C')+B'g_\q(A',C')-C'g_\q(A',B')\\ 
		&\quad +g_\q([A',B'],C')-g_\q([B',C'],A')-g_\q([A',C'],B')
	\end{align*}
	We pull this function back to $P$ and look for the function on $N$ which pulls back to the same expression. 
	
	The first three terms are easily handled, using the compatibility of the twist construction with tensor products and contractions. However, for the remaining terms we must use the formula for twists of commutators of vector fields \cite{Swa2010}:
	\begin{equation*}
		[A',B']\sim_{\mc H}[A,B]+f_\h^{-1}\omega_\h(A,B)Z
	\end{equation*}
	We work out one of the terms explicitly, using tildes to denote horizontal lifts of vector fields on $N$ and hats for horizontal lifts from $\bar N$:
	\begin{align*}
		\pi_{\bar N}^*(g_\q([A,B],C))&=(\pi^*_{\bar N}g_\q)\big(\widehat{[A',B']},\widehat{C'}\big)\\
		&=(\pi_N^*g_\h)(\widetilde{[A,B]}+\pi^*_N(f_\h^{-1}\omega_\h(A,B))\tilde Z,\tilde C)\\
		&=\pi_N^*(g_\h([A,B],C)+f_\h^{-1}\omega_\h(A,B)g_\h(Z,C))\\
	\end{align*} 
	This shows that
	\begin{align*}
		\pi_{\bar N}^*(2g_\q(D^\q_{A'}B',C'))&=\pi^*_N\Big(2g(D_A B,C)\\
		&\ \ +f_\h^{-1}(\omega_\h(A,B)g(Z,C)-\omega_\h(B,C)g(A,Z)-\omega_\h(A,C)g(B,Z))\Big)
	\end{align*}
	as claimed.
\end{myproof}

Combining the two lemmata, we see that, under the HK/QK correspondence, the Levi-Civita connection $D^\q$ of the quaternionic K\"ahler metric $g_\q$ is $\mc H$-related (in the sense of the previous lemma) to the connection $\nabla^S\coloneqq D+S^\h+S^\q$, where $D$ is the Levi-Civita connection of the hyper-K\"ahler metric $g$. Now we turn to the curvature:

\begin{lem}\label{lem:curvaturetwist}
	Let $(N,g)$ be a manifold with twist data $(Z,\omega_\h,f_\h)$ such that $g$ is invariant, and $D$ the Levi-Civita connection, with curvature tensor $R$. Let $\nabla=D+S$ be a $Z$-invariant connection on $N$ and $\nabla'$ the $\mc H$-related connection on the twist manifold, with curvature $R^{\nabla'}$. Then, for $Z$-invariant vector fields $A,B,C\sim_{\mc H}A',B',C'$, we have:
	\begin{equation*}
		R^{\nabla'}(A',B')C'\sim_{\mc H} R(A,B)C+T(A,B)C
	\end{equation*}
	where $T$ is the $\End(TM)$-valued two-form given by
	\begin{equation}\label{eq:Ttensor}
		T(A,B)C\coloneqq (D_AS)_B C-(D_BS)_A C+[S_A,S_B]C-\frac{1}{f_\h}\omega_\h(A,B)(D_CZ+S_ZC)
	\end{equation}
\end{lem}
\begin{myproof}
	Using the curvature formula $R(A,B)C =[\nabla_A,\nabla_B]C-\nabla_{[A,B]}C$, our definition of $\mc H$-relatedness for connections means that $[\nabla'_{A'},\nabla'_{B'}]C'\sim_{\mc H}[\nabla_A,\nabla_B]C$, while $\nabla'_{[A',B']}C'\sim_{\mc H} \nabla_{[A,B]}C+\frac{1}{f_\h}\omega_\h(A,B)\nabla_Z C$. Since $C$ is $Z$-invariant and $D$ is torsion-free, we may rewrite $\nabla_Z C=D_CZ+S_ZC$.
\end{myproof}

In the setting of the HK/QK correspondence, we use $S=S^\h+S^\q$, and will denote the curvature tensor of $\nabla^S$ by $R^S$.

\clearpage

\section{Main theorems}
\label{sec:maintheorems}

The formulae deduced in the previous section determine, in principle, the curvature of the quaternionic K\"ahler metric in terms of hyper-K\"ahler data. In their current form, however, they are not particularly useful because they are too complicated to allow for effective computations. The main complication arises from the repeated appearance of the tensor field $S$. Our goal in the following is to rewrite these expressions in a more usable form. 

\subsection{Formula for the Levi-Civita connection}
\label{sec:twistingconnections}

We first study the Levi-Civita connection under the HK/QK correspondence, aiming to simplify the expression for $S=S^\h+S^\q$ given by \eqref{eq:S^h} and \eqref{eq:S^q}.

\begin{thm}\label{thm:connection}
	Let $(N,g,I_k)$ be a (pseudo-)hyper-K\"ahler manifold with Levi-Civita connection $D$, and rotating circle action generated by the vector field $Z$. Denote its image under the HK/QK correspondence by $(\bar N,g_\q)$. Then the Levi-Civita connection $D^\q$ of $g_\q$ is $\mc H$-related to $\nabla^S=D+S$, where, in the notation of \Cref{sec:background}, we have
	\begin{equation}\label{eq:conncorrfinal}
		S_A B=\frac{1}{2}\sum_{\mu=0}^3 \bigg(\frac{1}{f_\h} g(I_\mu I_\h A,B) I_\mu Z
		- \frac{1}{f_Z} \big(\alpha_\mu(A) I_\mu I_1 B+\alpha_\mu(B) I_\mu I_1 A\big)\bigg)
	\end{equation}
\end{thm}
\begin{myproof}
	This is a long computation involving the identities derived in \Cref{sec:notation}. Regarding both $\iota_Zg_\h \otimes \omega_\h$ and $D g_\h$ as linear maps $\Gamma(TN^{\otimes 3})\to C^\infty(N)$, we may write
	\begin{equation*}
		2 f_Z^2g_\h((S^\h_A B+S^\q_AB,C)=\bigg(\frac{f_Z^2}{f_\h}\iota_Zg_\h\otimes  \omega_\h
		- f_Z^2 D g_\h\bigg)\Big(\!C\otimes A\otimes B\!-\!A\otimes B\otimes C\!-\!B\otimes A\otimes C\!\Big)
	\end{equation*}
	Recalling from \Cref{sec:notation} that $\frac{f_Z^2}{f_\h}g_\h |_{\H Z}=g|_{\H Z}$ and making use of the shorthand $g_\alpha=\sum_{\mu=0}^3(\alpha_\mu)^2$, we may rewrite the $(0,3)$-tensor in the first parenthesis as
	\begin{align*}
		&\alpha_0\otimes \omega_\h -f_Z^2 D \Big(\frac{1}{f_Z}g+\frac{1}{f_Z^2}g_\alpha\Big)\\
		&\quad=\alpha_0\otimes \omega_\h + \d f_Z\otimes  g+2\frac{\d f_Z}{f_Z}\otimes g_\alpha -D g_\alpha\\
		&\quad =\alpha_0\otimes \omega_1+\alpha_0\otimes \d \alpha_0 - \alpha_1\otimes \omega_0 - \frac{2}{f_Z} \alpha_1\otimes g_\alpha
		-D g_\alpha \\
		&\quad =\alpha_0\otimes \omega_1 - \alpha_1\otimes \omega_0 - \frac{2}{f_Z}\alpha_1\otimes g_\alpha
		+ 2\alpha_0\otimes \omega_0(DZ,\cdot)-D g_\alpha
	\end{align*}
	where the final step uses \eqref{eq:dalpha}. Our next step is to rework the final two terms, which feature derivatives of $Z$. We have
	\begin{align*}
		(D g_\alpha)(A\otimes B\otimes C)&=\sum_\mu (D_A\alpha_\mu)(B)\alpha_\mu (C)+\alpha_\mu (B)(D_A\alpha_\mu)(C)\\
		&=\bigg(\sum_\mu \alpha_\mu\otimes \omega_\mu(DZ,\cdot ) \bigg)\big(C\otimes A\otimes B+B\otimes A\otimes C\big)
	\end{align*}
	and thus, using \eqref{eq:dalpha} and the fact that $Z$ generates a rotating circle symmetry, we find
	\begin{align*}
		&\big(2\alpha_0\otimes \omega_0(DZ,\cdot)- Dg_\alpha\big)\big(C\otimes A\otimes B-A\otimes B\otimes C-B\otimes A\otimes C\big)\\
		&\quad =\big(\alpha_0\otimes \omega_0(DZ,\cdot)\big) \Big(2C\otimes A\otimes B-B\otimes C\otimes A-A\otimes C\otimes B\\
		&\hspace{3.8cm} -2A\otimes B\otimes C+B\otimes A\otimes C+C\otimes A\otimes B\\
		&\hspace{3.8cm} -2B\otimes A\otimes C+A\otimes B\otimes C+C\otimes B\otimes A\Big)\\
		&\qquad  -\bigg(\sum_{k=1}^3\alpha_k\otimes \omega_k(DZ,\cdot)\bigg)\Big(B\otimes C\otimes A+A\otimes C\otimes B-B\otimes A\otimes C\\
		&\hspace{5cm} -C\otimes A\otimes B-A\otimes B\otimes C-C\otimes B\otimes A\Big)\\
		&\quad =2\alpha_0(C)\omega_0(D_AZ,B) +2\alpha_1(C)\omega_1(D_AZ,B)\\
		&\qquad  +2\alpha_2(C)\omega_2(D_AZ,B) -\alpha_2(C)\omega_3(A,B)+\alpha_2(A)\omega_3(B,C)+\alpha_2(B)\omega_3(A,C) \\
		&\qquad  +2\alpha_3(C)\omega_3(D_AZ,B) +\alpha_3(C)\omega_2(A,B)-\alpha_3(A)\omega_2(B,C)-\alpha_3(B)\omega_2(A,C) \\
		&\quad =2\sum_\mu \alpha_\mu \otimes \omega_\mu(DZ,\cdot)(C\otimes A\otimes B)\\
		&\qquad \quad -(\alpha_2\otimes \omega_3-\alpha_3\otimes \omega_2)\big(C\otimes A\otimes B-A\otimes B\otimes C-B\otimes A\otimes C\big)
	\end{align*}
	With this, we arrive at
	\begin{align*}
		&2f_Z^2g_\h(S_AB,C)\\
		&\quad =2\bigg(\sum_\mu \alpha_\mu \otimes \omega_\mu (DZ,\cdot)\bigg)(C\otimes A\otimes B)\\
		&\qquad +\bigg(\alpha_0\otimes\omega_1-\alpha_1\otimes\omega_0-\alpha_2\otimes\omega_3+\alpha_3\otimes\omega_2
		-\frac{2}{f_Z}\alpha_1\otimes g_\alpha\bigg)\\
		&\hspace{7.5cm}\big(C\otimes A\otimes B - A\otimes B\otimes C - B\otimes A\otimes C\big)\\
		&\quad =2\bigg(\sum_\mu \alpha_\mu \otimes g(I_\mu DZ,\cdot)\bigg)(C\otimes A\otimes B)\\
		&\qquad +\bigg(\sum_\mu \alpha_\mu\otimes  g(I_\mu I_1\cdot,\cdot)-\frac{2}{f_Z}\alpha_1\otimes g_\alpha\bigg)\big(C\otimes A\otimes B-A\otimes B\otimes C-B\otimes A\otimes C\big)
	\end{align*}
	We rewrite this as
	\begingroup
	\allowdisplaybreaks
	\begin{align*}
		&\bigg(\sum_\mu \alpha_\mu \otimes g(I_\mu I_\h \cdot,\cdot)
		-\frac{2}{f_Z}\alpha_1\otimes g_\alpha\bigg)(C\otimes A\otimes B)\\
		&\quad -\bigg(\sum_\mu \alpha_\mu \otimes g(I_\mu I_1\cdot,\cdot)-\frac{2}{f_Z}\alpha_1\otimes g_\alpha\bigg)\big(A\otimes B\otimes C + B\otimes A\otimes C\big)\\
		&\quad =\sum_\mu\bigg( g(I_\mu I_\h A,B) g(I_\mu Z,C)-\frac{2}{f_Z}\alpha_\mu(A)\alpha_\mu(B) g(I_1Z,C)\\
		&\qquad -\alpha_\mu(A)g(I_\mu I_1B,C)-\alpha_\mu(B)g(I_\mu I_1 A,C)\\
		&\qquad +\frac{2}{f_Z}\big(\alpha_1(A)\alpha_\mu(B)+\alpha_1(B)\alpha_\mu(A)\big)g(I_\mu Z,C)\bigg)\\
		&\quad=g\bigg(\sum_\mu \bigg( g(I_\mu I_\h A,B)I_\mu Z-\frac{2}{f_Z}\alpha_\mu(A)\alpha_\mu(B) I_1 Z-\alpha_\mu(A)I_\mu I_1B-\alpha_\mu(B)I_\mu I_1 A\\
		&\qquad +\frac{2}{f_Z}\big(\alpha_1(A)\alpha_\mu(B)+\alpha_1(B)\alpha_\mu(A)\big)I_\mu Z\bigg) ,C\bigg)
	\end{align*}
	\endgroup
	Now, applying the endomorphism $\mc K$ to the first argument will bring this to the form $2f_Z^2g_\h(S^\h_AB+S^\q_AB,C)$, from which we can read off $S^\h_AB+S^\q_AB$. Using the formula \eqref{eq:metriccompareformula} for $\mc K$, we find
	\begin{align*}
		2f_Z^2 S_AB
		&=\frac{f_Z^2}{f_\h}\sum_\mu g(I_\mu I_\h A,B)I_\mu Z
		-f_Z\sum_\mu \big(\alpha_\mu(A)I_\mu I_1B + \alpha_\mu(B)I_\mu I_1 A\big)\\
		&\quad +\frac{2f_Z}{f_\h}\sum_\mu \big(\alpha_1(A)\alpha_\mu(B)I_\mu Z+\alpha_\mu(A)\alpha_1(B)I_\mu Z-\alpha_\mu(A)\alpha_\mu(B) I_1 Z\big)\\
		&\quad +\frac{f_Z}{f_\h} 
		\sum_{\mu,\lambda}\big(\alpha_\mu(A) \alpha_\lambda(I_\mu I_1B) + \alpha_\mu(B) \alpha_\lambda(I_\mu I_1 A)\big)I_\lambda Z
	\end{align*} 
	Dividing by $2f_Z^2$ and exchanging $\mu$ and $\lambda$ in the double summation, we obtain
	\begin{align*}
		S_AB
		&=\frac{1}{2f_\h}\sum_\mu g(I_\mu I_\h A,B)I_\mu Z
		-\frac{1}{2f_Z}\sum_\mu \big(\alpha_\mu(A)I_\mu I_1B + \alpha_\mu(B)I_\mu I_1 A\big)\\
		&\quad +\frac{1}{f_Z f_\h}\sum_\mu \big(\alpha_1(A)\alpha_\mu(B)I_\mu Z+\alpha_\mu(A)\alpha_1(B)I_\mu Z-\alpha_\mu(A)\alpha_\mu(B) I_1 Z\big)\\
		&\quad +\frac{1}{2f_Zf_\h} 
		\sum_{\mu,\lambda}\big(\alpha_\lambda (A) \alpha_\mu(I_\lambda I_1B) + \alpha_\lambda(B) \alpha_\mu(I_\lambda I_1 A)\big)I_\mu Z
	\end{align*}
	Let us study the last line in detail. We use the fact that we can replace $I_\lambda$ by $I_\lambda I_1$ in the summation over $\lambda$; this introduces no additional signs because $I_\lambda$ appears twice in the expression. Therefore we can rewrite it as
	\begin{align*}
		&-\frac{1}{2f_Zf_\h}\sum_{\mu,\lambda} \big(g(I_\lambda I_1Z,A)g(I_\mu Z,I_\lambda B)
		+g(I_\lambda I_1 Z,B)g(I_\mu Z,I_\lambda A)\big)I_\mu Z\\
		&\quad =-\frac{1}{2f_Zf_\h}\sum_{\mu,\lambda} \big(g(I_\lambda I_1Z,A)g(I_\lambda^{-1} I_\mu Z, B)
		+g(I_\lambda I_1 Z,B)g(I_\lambda^{-1} I_\mu Z, A)\big)I_\mu Z\\
		&\quad =\frac{1}{2f_Zf_\h}\sum_{\mu,\lambda} \big(g(I_\lambda I_1Z,A)g(I_\lambda I_\mu Z, B)
		+g(I_\lambda I_1 Z,B)g(I_\lambda I_\mu Z, A)\big)I_\mu Z\\
		&\qquad -\frac{1}{f_Zf_\h}\sum_{\mu}\big(\alpha_1(A)\alpha_\mu(B)
		+\alpha_1(B)\alpha_\mu(A)\big)I_\mu Z
	\end{align*}
	where in the last step we used the fact that $I^{-1}_k=-I_k$ for $k\in \{1,2,3\}$, while $I_0^{-1}=I_0$. The second line cancels out immediately in the expression for $S^\h_AB+S^\q_AB$, leaving us with
	\begin{align*}
		S_AB
		&=\frac{1}{2f_\h}\sum_\mu g(I_\mu I_\h A,B)I_\mu Z
		-\frac{1}{2f_Z}\sum_\mu \big(\alpha_\mu(A)I_\mu I_1B + \alpha_\mu(B)I_\mu I_1 A\big)\\
		&\quad -\frac{1}{f_Z f_\h}\sum_\mu \alpha_\mu(A)\alpha_\mu(B) I_1 Z\\
		&\quad +\frac{1}{2f_Zf_\h}\sum_{\mu,\lambda} 
		\big(g(I_\lambda I_1Z,A)g(I_\lambda I_\mu Z, B)
		+g(I_\lambda I_1 Z,B)g(I_\lambda I_\mu Z, A)\big)I_\mu Z\\
	\end{align*}
	Once again, we focus on the last line. First, we replace $I_\mu$ by $I_\mu I_1$ (as before, no additional signs appear), and consider the expression $\sum_\lambda g(I_\lambda I_1 Z,A)g(I_\lambda I_\mu I_1 Z,B)$, for fixed $\mu$. When $\mu=0$, this is manifestly symmetric in $(A,B)$. If $\mu\neq 0$, however, the substitution $I_\lambda\mapsto I_\lambda I_\mu$ yields
	\begin{equation*}
		\sum_\mu g(I_\lambda I_1 Z,A)g(I_\lambda I_\mu I_1 Z,B)=-\sum_\lambda g(I_\lambda I_\mu I_1 Z,A)g(I_\lambda I_1 Z,B)
	\end{equation*}
	which shows that the expression is antisymmetric in $(A,B)$. Therefore, all terms from the last line except those corresponding to $\mu=0$ vanish, but for $\mu=0$ they cancel in the expression for $S^\h_AB +S^\q_AB$. In summary, we have shown:
	\begin{equation*}
		S_AB=\frac{1}{2}\sum_\mu \bigg(\frac{1}{f_\h}g(I_\mu I_\h A,B)I_\mu Z-\frac{1}{f_Z}\big(\alpha_\mu(A)I_\mu I_1 B+\alpha_\mu(B)I_\mu I_1 A\big)\bigg)
	\end{equation*}
	This completes the proof.
\end{myproof} 

\begin{rem}
	Note that, in order to arrive at this expression, we have thoroughly mixed the contributions of $S^\h$ and $S^\q$. Thus, it seems that treating both the elementary deformation and the twist simultaneously is crucial. 
\end{rem}

\subsection{Curvature formula}
\label{sec:twistingcurvature}

As explained in \Cref{sec:twistformulae}, the curvature tensor $R^\q$ of the quaternionic K\"ahler metric $g_\q$ on $\bar N$ is $\mc H$-related to $\tilde R\coloneqq R^S-\frac{1}{f_\h}\omega_\h\otimes (DZ+S_Z)$, where $S$ is given by \Cref{thm:connection}, and $R^S$ is the curvature of $\nabla^S=D+S$, where $D$ is the Levi-Civita connection of the hyper-K\"ahler metric $g$ on $N$.

To state the final result, it is convenient to introduce the following pieces of notation:

\begin{mydef}\leavevmode
	\begin{numberedlist}
		\item We define the Kulkarni--Nomizu map
		\begin{equation*}
			\begin{tikzcd}[row sep=0]
				\Gamma\big((T^*N)^{\otimes 4}\big)\ar[r] & \Gamma\big(\bigwedge^2T^*N\otimes \bigwedge^2T^*N\big)\\
				\Phi \ar[r,mapsto] & \Phi^{\owedge}
			\end{tikzcd}
		\end{equation*}
		by setting 
		\begin{equation*}
			\Phi^\owedge(A,B,C,X)\coloneqq \Phi(A,C,B,X)-\Phi(A,X,B,C)+\Phi(B,X,A,C)-\Phi(B,C,A,X)
		\end{equation*}
		for arbitrary vector fields $A,B,C,X$.
		\item We define a second map
		\begin{equation*}
			\begin{tikzcd}[row sep=0]
				\Gamma\big(\bigwedge^2T^*N\otimes\bigwedge^2T^*N\big)\ar[r] & \Gamma\big(\bigwedge^2T^*N\otimes\bigwedge^2T^*N\big)\\
				\Phi \ar[r,mapsto] & \Phi^{\obar}
			\end{tikzcd}
		\end{equation*}
		by setting
		\begin{equation*}
			\Phi^\obar (A,B,C,X)\coloneqq \Phi^\owedge(A,B,C,X)+2\Phi(A,B,C,X)+2\Phi(C,X,A,B)
		\end{equation*}
	\end{numberedlist}
\end{mydef}

For $(0,2)$-tensors $\alpha$ and $\beta$, we set $\alpha\owedge\beta\coloneqq (\alpha\otimes\beta)^\owedge$ and analogously define $\alpha\obar\beta$. Note that our first definition then reduces to the usual Kulkarni--Nomizu product. 

It is well-known that the Kulkarni--Nomizu product of two symmetric bilinear forms is an algebraic curvature tensor, i.e.~possesses the algebraic symmetries of a (lowered) Riemann curvature tensor. When applied to two-forms, however, the resulting $(0,4)$-tensor generally fails to satisfy the Bianchi identity. This is remedied by $\obar$, which constructs an algebraic curvature tensor out of two two-forms.

Our main result is the following formula, which completely determines the curvature of the quaternionic K\"ahler metric:

\begin{thm}\label{thm:curvature}
	The Riemann curvature of the quaternionic K\"ahler metric $g_\q$ on $\bar N$ is $\mc H$-related to $\tilde R$, which is given by the expression
	\begin{align*}
		g_\h(\tilde R(A,B)C,X)&=\frac{1}{f_Z}g(R(A,B)C,X)\\
		&\quad + \frac{1}{8}\Big(g_\h\owedge g_\h+\sum_k g_\h(I_k\cdot,\cdot)\obar g_\h(I_k \cdot, \cdot)\Big) (A,B,C,X)\\
		&\quad-\frac{1}{8 f_Z f_\h}\Big(\omega_\h\obar \omega_\h+\sum_k \omega_\h(I_k \cdot,\cdot)\owedge\omega_\h(I_k \cdot,\cdot)\Big)(A,B,C,X)
	\end{align*}
	where $R$ is the curvature tensor of the initial (pseudo-)hyper-K\"ahler metric on $N$, and $A,B,C,X$ are arbitrary vector fields.
\end{thm}

\begin{rem}
	Each $I_k$, $k=1,2,3$, is skew with respect to $g_\h$, and the properties of $I_\h$ imply that $\omega_\h(I_k\cdot,\cdot)$ is a symmetric bilinear form. Thus, each individual term in the expression for $\tilde R$ is an algebraic curvature tensor.
\end{rem}

We divide the proof up into two steps. First, we give expressions for the various terms that make up $\tilde R$ (cf.~\eqref{eq:Ttensor}). We do this in the form of three lemmata. Then we sum them all up and further simplify to obtain the final curvature formula.

\begin{lem}\label{lem:termone}
	With $S$ as in \Cref{thm:connection} and $D$ the Levi-Civita connection of $g$, we have
	\begin{align*}
		&(D_AS)_BC-(D_BS)_AC\\
		&=\frac{1}{2}\sum_\mu \Bigg[ \frac{1}{f_\h^2}\big(\omega_\h(Z,A)\omega_\mu (I_\h B,C)I_\mu Z
		+\frac{1}{f_Z^2}\alpha_1(A)\big( g(I_\mu I_1 Z,B)I_\mu C+g(I_\mu I_1 Z,C)I_\mu B \big)\\
		&\quad \qquad \ \ -\frac{1}{f_\h^2}\omega_\h(Z,B)\omega_\mu (I_\h A,C)I_\mu Z
		-\frac{1}{f_Z^2}\alpha_1(B)\big( g(I_\mu I_1 Z,A)I_\mu C+g(I_\mu I_1 Z,C)I_\mu A \big)\\
		&\quad \qquad \ \ +\frac{1}{f_\h}\big(\omega_\mu(I_\h B,C)I_\mu D_AZ - \omega_\mu(I_\h A,C)I_\mu D_BZ
		+2\omega_\mu (R(A,B)Z,C)I_\mu Z\big)\\
		&\quad \qquad \ \ +\frac{1}{2f_Z}\Big(\big(g(I_\mu I_1 I_\h A,C)+g(I_\mu A,C)\big)I_\mu B
		-\big(g(I_\mu I_1I_\h B,C)+g(I_\mu B,C)\big)I_\mu A\Big)\\
		&\quad \qquad \ \ +\frac{1}{2f_Z}\big(\omega_\mu(A,B)-\omega_\mu(B,A)\big)I_\mu C\Bigg]-\frac{1}{2f_Z}\omega_\h(A,B) I_1 C
	\end{align*}
	for arbitrary vector fields $A,B,C$ on $N$.
\end{lem}
\begin{myproof}
	We first compute the covariant derivative of $S$, and anti-symmetrize later:
	\begin{align*}
		(D_A S)_B C&=\frac{1}{2}\sum_\mu \bigg(\!\!-\frac{\d f_\h(A)}{f_\h^2}g(I_\mu I_\h B,C)I_\mu Z+\frac{2}{f_\h}g(I_\mu (D_A(D_B Z)-D_{D_AB}Z),C)I_\mu Z\\
		&\qquad \qquad +\frac{1}{f_\h}g(I_\mu I_\h B,C)I_\mu D_A Z +\frac{\d f_Z(A)}{f_Z^2}(\alpha_\mu(B)I_\mu I_1 C+\alpha_\mu(C)I_\mu I_1 B)\\
		&\qquad \qquad -\frac{1}{f_Z}(g(I_\mu D_A Z,B)I_\mu I_1 C+g(I_\mu D_AZ,C)I_\mu I_1 B)\bigg)\\
		&=\frac{1}{2}\sum_\mu \bigg( \frac{1}{f_\h^2}\omega_\h(Z,A)\omega_\mu (I_\h B,C)I_\mu Z\\
		&\qquad \qquad +\frac{1}{f_Z^2}\alpha_1(A)(g(I_\mu I_1 Z,B)I_\mu C+g(I_\mu I_1Z,C)I_\mu B)\\
		&\qquad \qquad +\frac{1}{f_\h}(\omega_\mu(I_\h B,C)I_\mu D_AZ + 2 \omega_\mu (D_A(D_BZ)-D_{D_AB}Z,C))I_\mu Z\\
		&\qquad \qquad -\frac{1}{f_Z}\Big(\omega_\mu (D_A Z,B)I_\mu I_1 C-\frac{1}{2}g(I_\mu I_1(I_\h-I_1)A,C)I_\mu B\Big)\bigg)
	\end{align*}
	where, in the second step, we used the definitions of $f_Z$ and $f_\h$ as Hamiltonians, replaced $I_\mu$ by $I_\mu I_1$ in some of the terms, and used $DZ=\frac{1}{2}(I_\h-I_1)$ on the final term. We then grouped the terms according to their prefactors.
	
	After anti-symmetrizing in $(A,B)$, some of the terms involving covariant derivatives of $Z$ can be simplified considerably. Firstly, 
	\begin{equation*}
		D_A(D_BZ)-D_{D_AB}Z-D_B(D_AZ)+D_{D_BA}Z=R(A,B)Z
	\end{equation*}
	Secondly, we may apply the identity \eqref{eq:sumidentity} to the first term on the last line.
	
	Taking these two points into account, the claimed expression appears directly upon anti-symmetrization.
\end{myproof}

\begin{lem}\label{lem:termtwo}
	In the same notation as the previous lemma we have
	\begin{align*}
		[S_A,S_B]C&=\frac{1}{4}\sum_{\mu,\lambda}\bigg(
		\frac{1}{f_\h^2}\big(g(I_\mu I_\h A, Z)g(I_\lambda I_\h B,C)-g(I_\mu I_\h B, Z)g(I_\lambda I_\h A,C)\big)I_\lambda I_\mu Z\\
		&\hspace{1.5cm}+\frac{1}{f_Z^2}\Big[
		\big(g(I_\mu I_1 Z,A)g(I_\lambda I_1 Z,B)-g(I_\mu I_1 Z,B)g(I_\lambda I_1 Z,A)\big)I_\mu I_\lambda C\\
		&\hspace{1.5cm}+g(I_\mu I_1 Z,B)g(I_\lambda I_1 Z,C)I_\lambda I_\mu A 
		- g(I_\mu I_1 Z,A) g(I_\lambda I_1 Z,C)I_\lambda I_\mu B\Big]\bigg)\\
		&\quad+\frac{1}{4f_Zf_\h}\sum_\mu \Big( 2g(I_\h A,B)g(I_\mu I_1 Z,C)I_\mu Z\\
		&\hspace{2.7cm} -g(Z,Z)\big(g(I_\mu I_\h B,C)I_\mu I_1 A-g(I_\mu I_\h A,C)I_\mu I_1 B\big)\Big)
	\end{align*}
\end{lem}
\begin{myproof}
	Before getting started, we rewrite the formula for $S$ by replacing $I_\mu$ by $I_\mu I_1$ in the second part, so that we have 
	\begin{equation*}
		S_AB=\frac{1}{2}\sum_\mu \bigg(\frac{1}{f_\h}g(I_\mu I_\h A,B)I_\mu Z
		+\frac{1}{f_Z}\big(g(I_\mu I_1 Z,A)I_\mu B+g(I_\mu I_1 Z,B)I_\mu A\big)\bigg)
	\end{equation*}
	Now we apply this formula, replacing $B$ by $S_BC$:
	\begin{align*}
		S_A(S_B C)&=\frac{1}{2}\sum_\mu \bigg(\frac{1}{f_\h} g(I_\mu I_\h A,S_B C) I_\mu Z\\
		&\hspace{1.5cm}+ \frac{1}{f_Z} \big(g(I_\mu I_1Z,A) I_\mu S_B C+g(I_\mu I_1 Z, S_B C) I_\mu A\big)\bigg)\\
		&=\frac{1}{4}\sum_{\mu,\lambda}\bigg(
		\frac{1}{f_\h^2}g(I_\mu I_\h A,I_\lambda Z)g(I_\lambda I_\h B,C)I_\mu Z\\
		&\hspace{1.5cm}+ \! \frac{1}{f_Zf_\h}\Big[
		\big(g(I_\mu I_\h A,I_\lambda C)g(I_\lambda I_1 Z,B)
		\! + \! g(I_\mu I_\h A,I_\lambda B)g(I_\lambda I_1 Z,C) \big)I_\mu Z\\
		&\hspace{1.5cm}+g(I_\mu I_1 Z,A)g(I_\lambda I_\h B,C)I_\mu I_\lambda Z+g(I_\mu I_1 Z,I_\lambda Z)g(I_\lambda I_\h B,C)I_\mu A\Big]\\
		&\hspace{1.5cm} +\frac{1}{f_Z^2}\Big[g(I_\mu I_1 Z,A)\big( g(I_\lambda I_1 Z,B) I_\mu I_\lambda C + g(I_\lambda I_1 Z,C) I_\mu I_\lambda B \big)\\
		&\qquad \qquad \quad +\big(g(I_\mu I_1 Z,I_\lambda C)g(I_\lambda I_1 Z,B)+g(I_\mu I_1 Z,I_\lambda B)g(I_\lambda I_1 Z,C)\big)I_\mu A\Big]\bigg)
	\end{align*}
	We can rewrite all terms of the form $g(I_\mu X,I_\lambda Y)I_\mu W$ by replacing $I_\mu$ by $I_\lambda I_\mu$ as follows:
	\begin{align*}
		\sum_\mu g(I_\mu X, I_\lambda Y)I_\mu W&=\sum_\mu g(X,I_\mu^{-1}I_\lambda Y)I_\mu W
		=\sum_\mu g(X,I^{-1}_\mu Y)I_\lambda I_\mu W\\
		&=\sum_\mu g(I_\mu X,Y)I_\lambda I_\mu W
	\end{align*}
	Doing this several times, and also exchanging the labels $\mu$ and $\lambda$ in some terms, we obtain
	\begin{align*}
		S_A(S_BC)&=\frac{1}{4}\sum_{\mu,\lambda}\bigg(
		\frac{1}{f_\h^2}g(I_\mu I_\h A,Z)g(I_\lambda I_\h B,C)I_\lambda I_\mu Z\\
		&\qquad \qquad +\frac{1}{f_Z^2}\Big[g(I_\mu I_1 Z,A) g(I_\lambda I_1 Z,B) I_\mu I_\lambda C + g(I_\lambda I_1 Z,A)g(I_\mu I_1 Z,C) I_\lambda I_\mu B \\
		&\qquad \qquad +\big(g(I_\mu I_1 Z, C)g(I_\lambda I_1 Z,B)+g(I_\mu I_1 Z, B)g(I_\lambda I_1 Z,C)\big)I_\lambda I_\mu A\Big]\\
		&\qquad \qquad +\frac{1}{f_Zf_\h}\Big[
		\big(g(I_\mu I_\h A, C)g(I_\lambda I_1 Z,B)+g(I_\mu I_\h A, B)g(I_\lambda I_1 Z,C) \big)I_\lambda I_\mu Z\\
		&\qquad \qquad +g(I_\lambda I_1 Z,A)g(I_\mu I_\h B,C)I_\lambda I_\mu Z+g(I_\mu I_1 Z, Z)g(I_\lambda I_\h B,C)I_\lambda I_\mu A\Big]\bigg)
	\end{align*}
	We note that, among the terms multiplied by $f_Z^{-2}$, the second and third combine into an expression symmetric in $(A,B)$. The same is true for the first and third terms which come with a factor $(f_Z f_\h)^{-1}$, Antisymmetrizing, these cancel out and we find
	\begin{align*}
		[S_A,S_B]C&=\frac{1}{4}\sum_{\mu,\lambda}\bigg(
		\frac{1}{f_\h^2}\big(g(I_\mu I_\h A, Z)g(I_\lambda I_\h B,C)-g(I_\mu I_\h B, Z)g(I_\lambda I_\h A,C)\big)I_\lambda I_\mu Z\\
		&\hspace{1.5cm}+\frac{1}{f_Z^2}\Big[
		\big(g(I_\mu I_1 Z,A)g(I_\lambda I_1 Z,B)-g(I_\mu I_1 Z,B)g(I_\lambda I_1 Z,A)\big)I_\mu I_\lambda C\\
		&\hspace{1.5cm}+g(I_\mu I_1 Z,B)g(I_\lambda I_1 Z,C)I_\lambda I_\mu A 
		- g(I_\mu I_1 Z,A) g(I_\lambda I_1 Z,C)I_\lambda I_\mu B\Big]\\
		&\qquad \qquad +\frac{1}{f_Zf_\h}\Big[
		\big(g(I_\mu I_\h A,B)-g(I_\mu I_\h B,A)\big) g(I_\lambda I_1 Z,C)I_\lambda I_\mu Z\\
		&\qquad \qquad \qquad \quad +g(I_\mu I_1 Z,Z)\big(g(I_\lambda I_\h B,C)I_\lambda I_\mu A-g(I_\lambda I_\h A,C)I_\lambda I_\mu B\big)\Big]\bigg)
	\end{align*}
	The final two lines can be further simplified: Since $I_\h$ is skew, $g(I_\mu I_\h  A,B)$ is symmetric in $(A,B)$ unless $\mu=0$, in which case it is anti-symmetric. Therefore, only this terms survives. Furthermore, $g(I_\mu I_1 Z,Z)$ vanishes unless $\mu=1$. We are now left with
	\begin{align*}
		[S_A,S_B]C&=\frac{1}{4}\sum_{\mu,\lambda}\bigg(
		\frac{1}{f_\h^2}\big(g(I_\mu I_\h A, Z)g(I_\lambda I_\h B,C)-g(I_\mu I_\h B, Z)g(I_\lambda I_\h A,C)\big)I_\lambda I_\mu Z\\
		&\hspace{1.5cm}+\frac{1}{f_Z^2}\Big[
		\big(g(I_\mu I_1 Z,A)g(I_\lambda I_1 Z,B)-g(I_\mu I_1 Z,B)g(I_\lambda I_1 Z,A)\big)I_\mu I_\lambda C\\
		&\hspace{1.5cm}+g(I_\mu I_1 Z,B)g(I_\lambda I_1 Z,C)I_\lambda I_\mu A 
		- g(I_\mu I_1 Z,A) g(I_\lambda I_1 Z,C)I_\lambda I_\mu B\Big]\bigg)\\
		&\quad+\frac{1}{4f_Zf_\h}\sum_\mu \Big( 2g(I_\h A,B)g(I_\mu I_1 Z,C)I_\mu Z\\
		&\hspace{2.7cm} -g(Z,Z)\big(g(I_\mu I_\h B,C)I_\mu I_1 A-g(I_\mu I_\h A,C)I_\mu I_1 B\big)\Big)
	\end{align*}
	which is the expression we were after.
\end{myproof}

This determines $R^S$. Now, there is one more term to compute, namely $-\frac{1}{f_\h}\omega_\h\otimes (DZ+S_Z)$. Since we understand $\omega_\h$ well, we only need to study the second part:

\begin{lem}\label{lem:termthree}
	The $(1,1)$-tensor field $DZ+S_Z$ is given by the following expression:
	\begin{align*}
		D_AZ+S_ZA&=\frac{1}{2}\bigg(I_\h - \frac{f_\h}{f_Z}I_1\bigg)A
		+\frac{1}{2}\sum_\mu \bigg(\frac{1}{f_\h} g(I_\mu I_\h Z,A)
		+\frac{1}{f_Z}g(I_\mu I_1 Z,A) \bigg)I_\mu Z\\
	\end{align*}
\end{lem}
\begin{myproof}
	This is the result of a relatively short computation:
	\begin{align*}
		D_AZ+S_ZA&=\frac{1}{2}(I_\h-I_1)A\\
		&\quad +\frac{1}{2}\sum_\mu \bigg(\frac{1}{f_\h} g(I_\mu I_\h Z,A) I_\mu Z
		+\frac{1}{f_Z} \big(g(I_\mu I_1Z,Z) I_\mu A+g(I_\mu I_1 Z,A) I_\mu Z\big)\bigg)\\
		&=\frac{1}{2}(I_\h-I_1)A-\frac{g(Z,Z)}{2f_Z}I_1 A\\
		&\quad +\frac{1}{2}\sum_\mu \bigg(\frac{1}{f_\h} g(I_\mu I_\h Z,A)
		+\frac{1}{f_Z}g(I_\mu I_1 Z,A) \bigg)I_\mu Z\\
		&=\frac{1}{2}\bigg(I_\h - \frac{f_\h}{f_Z}I_1\bigg)A
		+\frac{1}{2}\sum_\mu \bigg(\frac{1}{f_\h} g(I_\mu I_\h Z,A)
		+\frac{1}{f_Z}g(I_\mu I_1 Z,A) \bigg)I_\mu Z\\
	\end{align*}
	In the second step, we once again used that $g(I_\mu I_1 Z,Z)$ vanishes unless $\mu=1$, while the final step used $f_\h=f_Z+g(Z,Z)$.
\end{myproof}

\begin{myproof}[Proof of \Cref{thm:curvature}]
	Now that we have computed all the individual pieces that make up $\tilde R$, we put them all together. We group terms according to their prefactors, and obtain the following expression:
	\begingroup
	\allowdisplaybreaks
	\begin{align*}
		&\tilde R(A,B)C-R(A,B)C\\
		&\quad=(D_AS)_BC-(D_BS)_A C+[S_A,S_B]C-\frac{1}{f_\h}\omega_\h(A,B)(D_CZ+S_ZC)\\
		&\quad =\frac{1}{4f_\h^2}\Bigg(2\sum_\mu
		\Big[\big(\omega_\h(Z,A)\omega_\mu (I_\h B,C) -\omega_\h(Z,B)\omega_\mu (I_\h A,C)\big)I_\mu Z\\
		&\hspace{2.7cm} -\omega_\h(A,B)g(I_\mu I_\h Z,C)I_\mu Z\Big]\\
		&\hspace{1.7cm} +\sum_{\mu,\lambda}\big(g(I_\mu I_\h A, Z)g(I_\lambda I_\h B,C)-g(I_\mu I_\h B, Z)g(I_\lambda I_\h A,C)\big)I_\lambda I_\mu Z\Bigg)\\
		&\qquad +\frac{1}{4f_Z^2}\Bigg(2\sum_\mu \Big[
		\alpha_1(A)\big( g(I_\mu I_1 Z,B)I_\mu C+g(I_\mu I_1 Z,C)I_\mu B \big)\\
		&\hspace{3cm} -\alpha_1(B)\big( g(I_\mu I_1 Z,A)I_\mu C+g(I_\mu I_1 Z,C)I_\mu A \big)\Big]\\
		&\hspace{1.7cm} +\sum_{\mu,\lambda}\Big[
		\big(g(I_\mu I_1 Z,A)g(I_\lambda I_1 Z,B)-g(I_\mu I_1 Z,B)g(I_\lambda I_1 Z,A)\big)I_\mu I_\lambda C\\
		&\hspace{3.2cm} +g(I_\mu I_1 Z,B)g(I_\lambda I_1 Z,C)I_\lambda I_\mu A 
		- g(I_\mu I_1 Z,A) g(I_\lambda I_1 Z,C)I_\lambda I_\mu B\Big]\Bigg)\\
		&\qquad +\frac{1}{4f_Zf_\h}\Bigg(\sum_\mu \Big[ 
		2g(I_\h A,B)g(I_\mu I_1 Z,C)I_\mu Z
		-2\omega_\h(A,B)g(I_\mu I_1 Z,C)I_\mu Z\\
		&\hspace{3.4cm} -g(Z,Z)\big(g(I_\mu I_\h B,C)I_\mu I_1 A
		-g(I_\mu I_\h A,C)I_\mu I_1 B\big)\Big]\Bigg)\\
		&\qquad +\frac{1}{2f_\h}\Bigg(
		\sum_\mu \Big[\omega_\mu(I_\h B,C)I_\mu D_AZ 
		- \omega_\mu(I_\h A,C)I_\mu D_BZ
		+2\omega_\mu (R(A,B)Z,C)I_\mu Z\Big]\\
		&\hspace{2.2cm}-\omega_\h(A,B)I_\h C\Bigg)\\
		&\qquad +\frac{1}{4f_Z}\Bigg(
		\sum_\mu \Big[\big(g(I_\mu I_1 I_\h A,C)+g(I_\mu A,C)\big)I_\mu B
		-\big(g(I_\mu I_1I_\h B,C)+g(I_\mu B,C)\big)I_\mu A\\
		&\hspace{3cm}+\big(\omega_\mu(A,B)-\omega_\mu(B,A)\big)I_\mu C\Big]\Bigg)
	\end{align*}
	\endgroup
	Note that the final terms of the expressions computed in \Cref{lem:termone} and \Cref{lem:termthree} canceled each other out.
	
	Now consider the terms proportional to $(f_Zf_\h)^{-1}$. The first two terms cancel, and we are left with
	\begin{align*}
		&-\frac{g(Z,Z)}{4f_Z f_\h}\sum_\mu \big( g(I_\mu I_\h B,C)I_\mu I_1 A -g(I_\mu I_\h A,C)I_\mu I_1 B \big)\\
		&\quad =\frac{1}{4}\bigg(\frac{1}{f_\h}-\frac{1}{f_Z}\bigg)\sum_\mu \big(\omega_\mu (I_\h B,C)I_\mu I_1 A
		-\omega_\mu (I_\h A,C)I_\mu I_1 B\big)\\
		&\quad =\frac{1}{4f_\h}\sum_\mu 
		\big(\omega_\mu (I_\h B,C)I_\mu (I_\h A-2D_AZ)-\omega_\mu (I_\h A,C)I_\mu (I_\h B-2D_BZ)\\
		&\qquad +\frac{1}{4f_Z}\sum_\mu 
		\big(g(I_\mu I_1I_\h B,C)I_\mu A-g (I_\mu I_1 I_\h A,C)I_\mu B\big)
	\end{align*}
	where, to get the last equality, we used $I_\h=I_1+2DZ$ and substituted $I_\mu I_1$ for $I_\mu$ in the last line. Now we compare to other terms from the expression for $\tilde R$. Observe that the entire second line is canceled out, as are the terms featuring $DZ$ on the first line.
	
	This leaves us with
	\begingroup
	\allowdisplaybreaks
	\begin{align*}
		&(D_AS)_BC-(D_BS)_A C+[S_A,S_B]C-\frac{1}{f_\h}\omega_\h(A,B)(D_CZ+S_ZC)\\
		&\quad =\frac{1}{4f_\h^2}\Bigg(2\sum_\mu
		\Big[\big(\omega_\h(Z,A)\omega_\mu (I_\h B,C) -\omega_\h(Z,B)\omega_\mu (I_\h A,C)\big)I_\mu Z\\
		&\hspace{2.9cm} -\omega_\h(A,B)g(I_\mu I_\h Z,C)I_\mu Z\Big]\\
		&\hspace{1.7cm} +\sum_{\mu,\lambda}\big(g(I_\mu I_\h A, Z)g(I_\lambda I_\h B,C)-g(I_\mu I_\h B, Z)g(I_\lambda I_\h A,C)\big)I_\lambda I_\mu Z\Bigg)\\
		&\qquad +\frac{1}{4f_Z^2}\Bigg(2\sum_\mu \Big[
		\alpha_1(A)\big( g(I_\mu I_1 Z,B)I_\mu C+g(I_\mu I_1 Z,C)I_\mu B \big)\\
		&\hspace{3.2cm} -\alpha_1(B)\big( g(I_\mu I_1 Z,A)I_\mu C+g(I_\mu I_1 Z,C)I_\mu A \big)\Big]\\
		&\hspace{1.7cm} +\sum_{\mu,\lambda}\Big[
		\big(g(I_\mu I_1 Z,A)g(I_\lambda I_1 Z,B)-g(I_\mu I_1 Z,B)g(I_\lambda I_1 Z,A)\big)I_\mu I_\lambda C\\
		&\hspace{2.9cm} +g(I_\mu I_1 Z,B)g(I_\lambda I_1 Z,C)I_\lambda I_\mu A 
		- g(I_\mu I_1 Z,A) g(I_\lambda I_1 Z,C)I_\lambda I_\mu B\Big]\Bigg)\\
		&\qquad +\frac{1}{4f_\h}\Bigg(
		\sum_\mu \Big[\omega_\mu(I_\h B,C)I_\mu I_\h A - \omega_\mu(I_\h A,C)I_\mu I_\h B
		+4\omega_\mu (R(A,B)Z,C)I_\mu Z\Big]\\
		&\hspace{2.2cm}-2\omega_\h(A,B)I_\h C\Bigg)\\
		&\qquad +\frac{1}{4f_Z}\Bigg(
		\sum_\mu \Big[g(I_\mu A,C)I_\mu B-g(I_\mu B,C)I_\mu A
		+\big(\omega_\mu(A,B)-\omega_\mu(B,A)\big)I_\mu C\Big]\Bigg)
	\end{align*}
	\endgroup
	Now we turn to the first group of terms, those multiplied by $f_\h^{-2}$. After splitting off the case $\lambda=0$ in the double summation, we can simplify:
	\begin{align*}
		&2\sum_\mu
		\Big[\omega_\h(Z,A)\omega_\mu (I_\h B,C) -\omega_\h(Z,B)\omega_\mu (I_\h A,C)
		-\omega_\h(A,B)\omega_\mu(I_\h Z,C)\Big]I_\mu Z\\
		&\quad +\sum_{\mu,\lambda}\big(g(I_\lambda I_\h A, Z)g(I_\mu I_\h B,C)
		-g(I_\lambda I_\h B, Z)g(I_\mu I_\h A,C)\big)I_\mu I_\lambda Z\\
		&\quad=\sum_\mu
		\Big[\omega_\h(Z,A)\omega_\mu (I_\h B,C) -\omega_\h(Z,B)\omega_\mu (I_\h A,C)
		-2\omega_\h(A,B)\omega_\mu(I_\h Z,C)\Big]I_\mu Z\\
		&\qquad +\sum_\mu\sum_{k=1}^3\big(g(I_k I_\h Z, A)\omega_\mu(I_\h B,C)
		-g(I_k I_\h Z, B)\omega_\mu (I_\h A,C)\big)I_\mu I_k Z\\
		&\quad =-2\sum_\mu \omega_\h (A,B)\omega_\mu(I_\h Z,C) I_\mu Z\\
		&\qquad +\sum_{\mu,\lambda}\big(\omega_\lambda(I_\h Z, A)\omega_\mu(I_\h B,C)
		-\omega_\lambda(I_\h Z,B)\omega_\mu(I_\h A,C)\big)I_\mu I_\lambda Z
	\end{align*}
	The $f_Z^{-2}$-terms may be simplified analogously. This yields:
	\begingroup
	\allowdisplaybreaks
	\begin{align*}
		&2\sum_\mu \Big[
		\alpha_1(A)g(I_\mu I_1 Z,B)-\alpha_1(B)g(I_\mu I_1 Z,A)\big)I_\mu C\\
		&\hspace{1cm} +\alpha_1(A)g(I_\mu I_1 Z,C)I_\mu B-\alpha_1(B)g(I_\mu I_1 Z,C)I_\mu A \big)\Big]\\
		&\quad +\sum_{\mu,\lambda}\Big[
		\big(g(I_\lambda I_1 Z,A)g(I_\mu I_1 Z,B)-g(I_\lambda I_1 Z,B)g(I_\mu I_1 Z,A)\big)I_\lambda I_\mu C\\
		&\hspace{1.5cm}+g(I_\lambda I_1 Z,B)g(I_\mu I_1 Z,C)I_\mu I_\lambda A 
		- g(I_\lambda I_1 Z,A) g(I_\mu I_1 Z,C)I_\mu I_\lambda B\Big]\\
		&\quad=\sum_{\mu,\lambda}\Big[
		\big(g(I_\mu I_1 Z,B)g(I_1 Z,I_\lambda A) - g(I_\mu I_1 Z,A)g(I_1 Z,I_\lambda B)\big)I_\mu I_\lambda C\\
		&\hspace{1.5cm} +g(I_1 Z,I_\lambda A) g(I_\mu I_1 Z,C)I_\mu I_\lambda B - g(I_1 Z,I_\lambda B)g(I_\mu I_1 Z,C)I_\mu I_\lambda A \Big]
	\end{align*}
	\endgroup
	Now, making the substitutions $I_\lambda \mapsto I_\lambda^{-1}$ and $I_\mu\mapsto I_\mu I_\lambda$, we obtain
	\begin{align*}
		&\sum_{\mu,\lambda}\Big[g(I_\mu I_\lambda I_1 Z,B)g(I_1 Z,I_\lambda^{-1}A) - g(I_\mu I_\lambda I_1 Z,A)g(I_1 Z,I_\lambda^{-1} B)\big)I_\mu C\\*
		&\qquad+g(I_1 Z,I_\lambda^{-1} A) g(I_\mu I_\lambda I_1 Z,C)I_\mu B - g(I_1 Z,I_\lambda^{-1} B)g(I_\mu I_\lambda I_1 Z,C)I_\mu A \Big]\\
		&\quad =\sum_{\mu,\lambda} \Big[g(I_\lambda I_1 Z,I_\mu^{-1} B)g(I_\lambda I_1 Z,A) - g(I_\lambda I_1 Z,I_\mu^{-1} A)g(I_\lambda I_1 Z, B)\big)I_\mu C\\*
		&\hspace{1.6cm}+g(I_\lambda I_1 Z, A) g(I_\lambda I_1 Z,I_\mu^{-1} C)I_\mu B - g(I_\lambda I_1 Z, B)g( I_\lambda I_1 Z,I_\mu^{-1} C)I_\mu A \Big]\\
		&=\sum_\mu \Big[\big(g_\alpha(I_\mu^{-1}B,A)-g_\alpha(I_\mu^{-1}A,B)\big)I_\mu C 
		+ g_\alpha(A,I_\mu^{-1}C)I_\mu B-g_\alpha(B,I_\mu^{-1}C)I_\mu A\Big]
	\end{align*}
	where the last step consists of replacing $I_\lambda$ by $I_\lambda I_1$ and recalling $g_\alpha\coloneqq \sum \alpha_\mu^2$.
	
	With these simplifications, we have arrived at
	\begingroup
	\allowdisplaybreaks
	\begin{align*}
		&\tilde R(A,B)C-R(A,B)C\\
		&\quad =\frac{1}{4f_\h^2}\bigg(\sum_{\mu,\lambda}\big(\omega_\lambda(I_\h Z, A)\omega_\mu(I_\h B,C)
		-\omega_\lambda(I_\h Z,B)\omega_\mu(I_\h A,C)\big)I_\mu I_\lambda Z\\
		&\hspace{1.8cm}-2\sum_\mu \omega_\h(A,B)\omega_\mu(I_\h Z,C)I_\mu Z\bigg)\\
		&\qquad +\frac{1}{4f_\h}\bigg(
		\sum_\mu \Big[\omega_\mu(I_\h B,C)I_\mu I_\h A
		-\omega_\mu(I_\h A,C)I_\mu I_\h B
		+4\omega_\mu (R(A,B)Z,C)I_\mu Z\Big]\\
		&\hspace{2.1cm} -2\omega_\h(A,B)I_\h C\bigg)\\
		&\qquad +\frac{1}{4f_Z^2}\sum_\mu \Big[\big(g_\alpha(B,I_\mu A)-g_\alpha(A,I_\mu B)\big)I_\mu C 
		+ g_\alpha(I_\mu A,C)I_\mu B-g_\alpha(I_\mu B,C)I_\mu A\Big]\\
		&\qquad +\frac{1}{4f_Z}\sum_\mu \Big[g(I_\mu A,C)I_\mu B-g(I_\mu B,C)I_\mu A
		+\big(\omega_\mu(A,B)-\omega_\mu(B,A)\big)I_\mu C\Big]\\
		&\quad =\frac{1}{4f_\h^2}\Bigg(\sum_{\mu,\lambda}\big(\omega_\mu(I_\h A,C)\alpha_\lambda(I_\h B)
		-\omega_\mu(I_\h B,C)\alpha_\lambda(I_\h A)\big)I_\mu I_\lambda Z\\
		&\qquad \qquad +2\sum_\mu \omega_\h(A,B)\alpha_\mu(I_\h C)I_\mu Z\Bigg)\\
		&\qquad +\frac{1}{4f_\h}\Bigg(
		\sum_\mu \Big[\omega_\mu(I_\h B,C)I_\mu I_\h A
		-\omega_\mu(I_\h A,C)I_\mu I_\h B
		-4\alpha_\mu (R(A,B)C)I_\mu Z\Big]\\
		&\qquad \qquad \qquad -2\omega_\h(A,B)I_\h C\Bigg)\\
		&\qquad +\frac{1}{4}\sum_\mu \Big[g_\h(I_\mu A,C)I_\mu B-g_\h(I_\mu B,C)I_\mu A+\big(g_\h(I_\mu A,B)-g_\h (I_\mu B,A)\big)I_\mu C\Big]
	\end{align*}
	\endgroup
	where we used that $R(A,B)$ and $I_\h$ both commute with each $I_\mu$ in the last step.
	
	Recalling \eqref{eq:metriccompareformula}, we have
	\begin{equation*}
		\frac{1}{f_Z}\mc K(A)=A -\frac{1}{f_\h}\sum_\mu \alpha_\mu(A)I_\mu Z
	\end{equation*}
	Applying this to the vector fields $R(A,B)C$, $I_\h A$ and $I_\h B$ and using that $I_\mu$ and $\mc K$ commute (cf.~\Cref{lem:metriccomparecommutes}), we find
	\begin{align*}
		&\tilde R(A,B)C\\
		&\quad =\frac{1}{f_Z}\mc K(R(A,B)C)-\frac{1}{2f_Zf_\h}\omega_\h(A,B)\mc K(I_\h C)\\\numberthis\label{eq:finalcurvexpression}
		&\qquad +\frac{1}{4f_Zf_\h}\sum_\mu\big(\omega_\mu(I_\h B,C) \mc K(I_\h I_\mu A)
		-\omega_\mu(I_\h A,C)\mc K(I_\h I_\mu B)\big)\\
		&\qquad +\frac{1}{4}\sum_\mu \Big[g_\h(I_\mu A,C)I_\mu B-g_\h(I_\mu B,C)I_\mu A+\big(g_\h(I_\mu A,B)-g_\h (I_\mu B,A)\big)I_\mu C\Big]
	\end{align*}
	To obtain the expression given in the statement of the theorem, we contract with an auxiliary vector field $X$, and split off the case $\mu=0$ in the summations:
	\begin{align*}
		&g_\h(\tilde R(A,B)C,X)\\
		&\quad =\frac{1}{f_Z}g(R(A,B)C,X)\\
		&\qquad +\frac{1}{4 f_Zf_\h}\Big[\omega_\h(B,C)\omega_\h(A,X)-\omega_\h(A,C)\omega_\h(B,X)-2\omega_\h(A,B)\omega_\h(C,X)\Big]\\
		&\qquad +\frac{1}{4f_Z f_\h}\sum_{k=1}^3 \Big[\omega_\h(I_k B,C)\omega_\h(I_k A,X)-\omega_\h(I_k A,C)\omega_\h (I_k B,X)\Big]\\
		&\qquad +\frac{1}{4}\big(g_\h(A,C)g_\h(B,X)-g_\h(B,C)g_\h(A,X)\big)\\
		&\qquad +\frac{1}{4}\sum_{k=1}^3\! \Big[g_\h(I_k A,C)g_\h(I_k B,X)\! - \! g_\h(I_k B,C)g_\h(I_k A,X)
		\! + \! 2g_\h(I_k A,B)g_\h(I_k C,X)\Big]\\
		&\quad =\frac{1}{f_Z}g(R(A,B)C,X)\\
		&\qquad -\frac{1}{8 f_Z f_\h}\Big(\omega_\h\obar \omega_\h+\sum_k \omega_\h(I_k \cdot,\cdot)\owedge\omega_\h(I_k \cdot,\cdot)\Big)  (A,B,C,X)\\
		&\qquad + \frac{1}{8}\Big(g_\h\owedge g_\h+\sum_k g_\h(I_k\cdot,\cdot)\obar g_\h(I_k \cdot, \cdot)\Big) (A,B,C,X)
	\end{align*}
	In the final step, we grouped terms into algebraic curvature tensors. This finishes the proof.
\end{myproof}

A well-known theorem of Alekseevsky~\cite{Ale1968} asserts that the curvature tensor of any quaternionic K\"ahler manifold $\bar N$ is of the form $R^\q=\nu R_0+R_1$, where $\nu=\frac{\scal}{4n(n+2)}$ is the reduced scalar curvature ($n$ being the quaternionic dimension), $R_0$ is the curvature tensor of quaternionic projective space $\HP^n$ and $R_1$ (or rather its lowered version) is an algebraic curvature tensor of hyper-K\"ahler type. The latter condition means that, for any two vector fields $X$ and $Y$, the endomorphism $R_1(X,Y)$ commutes with any section of the rank three bundle $Q\subset \End T\bar N$ defining the quaternionic K\"ahler structure. The tensor $R_1$ is sometimes called the quaternionic Weyl curvature.

With respect to its symmetric metric $g_0$ of reduced scalar curvature $1$, the curvature of $\HP^n$ is locally given by the expression
\begin{equation*}
	g_0(R_0(A,B)C,X)=-\frac{1}{8}\bigg(g_0\owedge g_0 
	+\sum_{k=1}^3 g_0(J_k\cdot,\cdot)\obar g_0(J_k\cdot,\cdot)\bigg)(A,B,C,X)
\end{equation*}
where $\{J_1,J_2,J_3\}$ is a local orthonormal frame for $Q$. Comparing to \Cref{thm:curvature}, this clearly corresponds to the second line of our curvature formula (note that $\nu=-1$). We can check explicitly that the remaining terms are of hyper-K\"ahler type.

\begin{lem}
	The algebraic curvature tensor $R_1$ on $\bar N$, determined by the relation
	\begin{align*}
		g_\q(R_1(A,B)C,X)&\sim_{\mc H}\frac{1}{f_Z}g(R(A,B)C,X)\\
		&\qquad-\frac{1}{8 f_Z f_\h}\Big(\omega_\h\obar \omega_\h+\sum_k \omega_\h(I_k \cdot,\cdot)\owedge\omega_\h(I_k \cdot,\cdot)\Big)  (A,B,C,X)
	\end{align*}
	is of hyper-K\"ahler type.
\end{lem}
\begin{myproof}
	We need to check that $R_1$ commutes with any section of $Q$. Since this is a pointwise condition, it suffices to check the assertion on a frame, which we may transfer to the hyper-K\"ahler side. There, we can work with respect to the convenient frame provided by $\{I_1,I_2,I_3\}$. Thus, it suffices to prove that the $\mc H$-related tensor field on $N$ commutes with each $I_k$. This is obviously the case for the first term, which is just the hyper-K\"ahler curvature tensor of $g$ (up to a pointwise scaling). For the remaining terms, we have:
	\begin{align*}
		&\Big(\omega_\h\obar \omega_\h+\sum_k \omega_\h(I_k \cdot,\cdot)\owedge\omega_\h(I_k \cdot,\cdot)\Big)  (A,B,I_j C,I_j X)\\
		&\quad =4\omega_\h(A,B)\omega_\h(I_j C, I_j X)\\
		&\qquad +2\sum_{\mu=0}^3 \Big(\omega_\h(I_\mu A,I_j C)\omega_\h(I_\mu B,I_j X)
		-\omega_\h(I_\mu A,I_j X)\omega_\h(I_\mu B,I_j C)\Big)\\
		&\quad =4\omega_\h(A,B)\omega_\h(C,X)\\
		&\qquad +2\sum_{\mu=0}^3 \Big(\omega_\h(I_j I_\mu A,C)\omega_\h(I_j I_\mu B, X)
		-\omega_\h(I_j I_\mu A,X)\omega_\h(I_j I_\mu B,C)\Big)\\
		&\quad =\Big(\omega_\h\obar \omega_\h+\sum_k \omega_\h(I_k \cdot,\cdot)\owedge\omega_\h(I_k \cdot,\cdot)\Big)  (A,B,C,X)\\
	\end{align*}
	where, in the final step, we replaced each instance of $I_\mu$ by $I_j^{-1}I_\mu$ in the summation. This suffices to prove the claim, since for any endomorphism $E$ and vector fields $X,I_kY$, we have $g([E,I_k]X,I_kY)=g(E I_k X,I_k Y)-g(EX,Y)$.
\end{myproof}

Thus, our results refine Alekseevsky's decomposition by providing a precise expression for the quaternionic Weyl curvature.

\clearpage

\section{Application: Curvature norm}
\label{sec:applications}

The main use of the curvature formula is to enable us to compute curvature quantities for quaternionic K\"ahler manifold that arise via the HK/QK correspondence. This approach is most effective if the curvature of the corresponding hyper-K\"ahler manifold is as simple as possible. Here, we study the extreme case where its curvature tensor vanishes identically.

For any non-negative integer $n$, we consider the manifold
\begin{equation*}
	N_n=\bigg\{(z_0,\dots,z_n,w_0,\dots,w_n)\in \C^{n+1}\times \C^{n+1}\bigg| -\abs{z_0}^2+\sum_{i=1}^n \abs{z_i}^2<c\bigg\}
\end{equation*}
where $c\geq 0$ is a real constant, and equip it with the standard flat (pseudo-)hyper-K\"ahler structure induced by
\begin{equation}\label{eq:flatHKstr}
	\begin{aligned}
		g&=-\big(\abs{\d z_0}^2+\abs{\d w_0}^2\big)+\sum_{j=1}^n \big(\abs{\d z_j}^2+\abs{\d w_j}^2\big)\\
		\omega_1&=\frac{i}{2}\bigg(-\big(\d z_0\wedge \d \bar z_0 + \d w_0 \wedge \d \bar w_0\big)
		+\sum_{j=1}^n \big(\d z_j \wedge \d \bar z_j + \d w_j \wedge \d \bar w_j\big)\bigg)\\
		\omega_2&=\frac{i}{2}\bigg(\d z_0 \wedge \d w_0 - \d \bar z_0 \wedge \d \bar w_0 
		+\sum_{j=1}^n \big(\d z_j\wedge \d w_j - \d \bar z_j \wedge \d \bar w_j\big)\bigg)\\
		\omega_3&=\frac{1}{2}\bigg(\d z_0 \wedge \d w_0 + \d \bar z_0 \wedge \d \bar w_0
		+\sum_{j=1}^n \big(\d z_j \wedge \d w_j + \d \bar z_j \wedge \d \bar w_j\big)\bigg)
	\end{aligned}
\end{equation}
The (inverse of the) standard action of $\C^*$ induces a rotating circle symmetry, generated by the vector field
\begin{equation*}
	Z=-i\sum_{j=0}^n \bigg(z_j \pd{}{z_j} - \bar z_j \pd{}{\bar z_j}\bigg)
\end{equation*}
It preserves $\omega_1$, with Hamiltonian function
\begin{equation*}
	f_Z=\frac{1}{2}\bigg(\abs{z_0}^2-\sum_{j=1}^n\abs{z_j}^2\bigg)-\frac{1}{2}c
\end{equation*}
which is everywhere positive for $c\geq 0$. 

To apply the HK/QK correspondence, we consider the elementary deformation $g_\h$ of the metric and the twist data $(Z,\omega_\h,f_\h)$, where
\begin{equation}\label{eq:flatomegaH}
	\omega_\h=\frac{i}{2}\bigg(\d z_0 \wedge \d \bar z_0 - \d w _0 \wedge \d \bar w_0 
	+\sum_{j=1}^n\big(-\d z_j \wedge \d \bar z_j + \d w_j \wedge \d \bar w_j\big)\bigg)
\end{equation}
and
\begin{equation*}
	f_\h=-\frac{1}{2}\bigg(\abs{z_0}^2-\sum_{i=j}^n\abs{z_j}^2\bigg)-\frac{1}{2}c
\end{equation*}

\begin{lem}
	On every $N_n$, the endomorphism fields $\mc K$, $I_\mu$, $\mu=0,1,2,3$, and $I_\h$ all commute, and $I_\h^2=-\id$.
\end{lem}
\begin{myproof}
	We already know that $\mc K$ and $I_\h$ each commute with every $I_\mu$, so all that remains is to prove that $I_\h$ commutes with $\mc K$. This will follow once we prove that $I_\h$ preserves the subbundle $\H Z$ of $TN$, as well as its orthogonal complement $(\H Z)^\perp$, since $\mc K$ restricts to multiples of the identity on these complementary subbundles. 
	
	In fact, it is enough to check that $g(I_\h Z,X)=0$ for any $X\in\Gamma((\H Z)^\perp)\subset \mf X(N)$. This follows from the fact that $I_\h$ is skew with respect to $g$ and commutes with each $I_\mu$, and that each $I_\mu$ acts $g$-orthogonally. We find
	\begin{equation*}
		g(I_\h Z,X)=\omega_\h(Z,X)=g(I_1 Z,X)+\d(\iota_Z g)(Z,X)
		=-X(g(Z,Z))
	\end{equation*}
	where we used the fact that $Z$ is a Killing field. Using the explicit coordinate expressions given for $g$ and $Z$, we find $g(Z,Z)=-\Big(\abs{z_0}^2-\sum_{j=1}^n\abs{z_j}^2\Big)$ and consequently
	\begin{equation*}
		\d(g(Z,Z))=-z_0\d \bar z_0-\bar z_0\d z_0 + \sum_{j=1}^n (z_j \d \bar z_j + \bar z_j \d z_j)
		=2\iota_{I_1 Z}g
	\end{equation*}
	This shows that $X(g(Z,Z))=0$, proving our claim. Thus, $\mc K$, $I_\h$ and $I_k$ all commute in this case. The fact that $I_\h^2=-\id$ can be seen by comparing the coordinate expressions \eqref{eq:flatHKstr} and \eqref{eq:flatomegaH} for $\omega_1$ and $\omega_\h$.
\end{myproof}

It is known that, for $c=0$, the HK/QK correspondence maps this flat hyper-K\"ahler manifold to the non-compact symmetric space $\bar N_n=\frac{SU(n+1,2)}{S(U(n+1)\times U(2))}$. For $c>0$, the resulting quaternionic K\"ahler metric is known to be complete~\cite{CDS2017}, but its isometry group is not well-understood. In \cite{CST2020}, we use general arguments to prove that it acts by cohomogeneity at most one. In the following, we will show that this result is sharp, i.e.~that the metric is of cohomogeneity exactly one. 

Concretely, we will use our curvature formula to compute the norm of the curvature, regarded as a self-adjoint endormorphism $\mc R^\q:\bigwedge^2 T\bar N_n\to \bigwedge^2 T\bar N_n$ (up to a factor of four, this is the same as the so-called Kretschmann scalar). With respect to an orthonormal frame $\{e_a\}$ of $T\bar N_n$, this curvature invariant can be expressed as follows:
\begin{equation*}
	\norm{\mc R^\q}^2=\sum g_\q(\mc R^\q(e_a \wedge e_b),e_c\wedge e_d)g_\q(\mc R^\q(e_c\wedge e_d),e_a\wedge e_b)
\end{equation*}
where the summation runs over the corresponding orthonormal basis of $\bigwedge^2T\bar N_n$.

Due the compatibility of the twist construction with respect to tensor products and contractions, we can compute this quantity on the hyper-K\"ahler side of the correspondence:
\begin{equation*}
	\norm{\mc R^\q}^2=\sum g_\h(\tilde{\mc R}(e'_a \wedge e'_b),e'_c\wedge e'_d)g_\h(\tilde{\mc R}(e'_c\wedge e'_d),e'_a\wedge e'_b)
\end{equation*}
where $\{e'_a\}$ is now an orthonormal frame with respect to $g_\h$. Our first step is to write this function invariantly. Using the metric $g_\h$, we can extend the definition of the Kulkarni--Nomizu product $\owedge$ to endomorphisms of $TN_n$. Indeed, for $E,F\in \End TN_n$ and vector fields $A,B,C,X$, we may define $E\owedge_{g_\h} F:\bigwedge^2 TN_n\to \bigwedge^2TN_n$ via
\begin{equation*}
	g_\h((E\owedge_{g_\h}F)(A\wedge  B),C\wedge X)\coloneqq \big(g_\h(E\cdot,\cdot)\owedge g_\h(F\cdot,\cdot)\big)(A,B,C,X)
\end{equation*}
Of course, we are mostly interested in the case where $E$ and $F$ are self-adjoint with respect to $g_\h$ (so that we obtain an algebraic curvature tensor). For skew-adjoint endomorphisms $E,F$, we may analogously define $E\obar_{g_\h}F:\bigwedge^2TN_n\to \bigwedge^2TN_n$. This notation allows us to succinctly express the operator $\tilde{\mc R}:\bigwedge^2TN_n\to \bigwedge^2TN_n$ associated to $\tilde R$:
\begin{align*}
	&g_\h(\tilde{\mc R}(A\wedge B),C\wedge X)\\
	&\hspace{1cm}=\frac{1}{f_Z}g(\mc R(A\wedge B),C\wedge X)\\
	& \hspace{1cm}\quad + \frac{1}{8}
	g_\h\bigg(\Big(\id\owedge_{g_\h} \id +\sum_k I_k\obar_{g_\h} I_k\Big)(A\wedge B),C\wedge X\bigg)\\
	&\hspace{1cm}\quad -\frac{1}{8 f_Z f_\h}g_\h\bigg(\!\Big(\mc KI_\h\obar_{g_\h}\mc K I_\h+\sum_k \mc K I_\h I_k \owedge_{g_\h}\mc K I_\h I_k\Big)(A\wedge B),C\wedge X\!\bigg)
\end{align*}

\begin{thm}
	Let $\mc R^\q$ denote the curvature endomorphism of the quaternionic K\"ahler manifold $(\bar N_{n-1},g_\q)$, $n\geq 1$, obtained by twisting $(N_{n-1},g_\h)$ with respect to the twist data $(Z,\omega_\h,f_\h)$. Then 
	\begin{equation*}
		\norm{\mc R^\q}^2=n(5n+1)
		+3\bigg(\frac{f_Z^3}{f_\h^3}+(n-1)\frac{f_Z}{f_\h}\bigg)^2
		+3\bigg(\frac{f_Z^6}{f_\h^6}+(n-1)\frac{f_Z^2}{f_\h^2}\bigg)
	\end{equation*}
\end{thm}
\begin{myproof}
	In this family of examples, the curvature $\mc R$ of the hyper-K\"ahler metric vanishes. We can therefore express the norm of $\mc R^\q$ as follows:
	\begin{equation*}
		\norm{\mc R^\q}^2\!=\!\frac{1}{64}\!\tr\! \Bigg[\!\bigg(\!\!
		\id\owedge_{g_\h} \id +\sum_k I_k\obar_{g_\h}I_k
		-\frac{1}{f_Z f_\h}\Big(\mc K I_\h\obar_{g_\h}\mc K I_\h+\sum_k \mc K I_\h I_k \owedge_{g_\h}\mc K I_\h I_k\Big)\!\bigg)^2\Bigg]
	\end{equation*}
	Thus, we have to compute three types of terms. The following lemma reduces this to the computation of traces of various compositions of the endomorphisms involved.
	
	\begin{lem}
		Let $h$ be a (pseudo-)Riemannian metric on a manifold $N$, let $E,F$ be self-adjoint endomorphisms of $TN$ and let $K,L$ be skew-adjoint endomorphisms of $TN$, with respect to $h$. Then
		\begin{numberedlist}
			\item $\tr\big((E\owedge_h E)\circ (F\owedge_h F)\big)
			=2\big(\tr (E\circ F)\big)^2-2\tr\big((E\circ F)^2\big)$.
			\item $\tr\big((K\obar_h K)\circ (L\obar_h L)\big)
			=6\big(\tr(K\circ L)\big)^2+6\tr\big((K\circ L)^2\big)$.
			\item $\tr\big((E\owedge_h E)\circ (K\obar_h K)\big)
			=2\big(\tr(E\circ K)\big)^2-6\tr\big((E\circ K)^2\big)$.
		\end{numberedlist}
	\end{lem}
	\begin{myproof}
		Let $\{e_a\}$ be an orthonormal basis for $h$ with $h(e_a,e_a)=\epsilon_a\in \{\pm 1\}$.
		\begin{numberedlist}
			\item Directly from the definition of $\owedge_h$, we have
			\begin{align*}
				&\tr\big(E\owedge_h E)\circ (F\owedge_h F)\big)\\
				&\quad =\frac{1}{4}\sum_{a,b,c,d}\epsilon_a\epsilon_b\epsilon_c\epsilon_d
				h((E\owedge_h E)e_c\wedge e_d,e_a\wedge e_b)
				h((F\owedge_h F)e_a\wedge e_b,e_c\wedge e_d)\\
				&\quad =\sum_{a,b,c,d}\epsilon_a\epsilon_b\epsilon_c\epsilon_d
				\big(h(E e_c,e_a)h(E e_d,e_b)-h(E e_c,e_b)h(E e_d,e_a)\big)\\
				&\hspace{3.1cm}\big(h(F e_a,e_c)h(F e_b,e_d)-h(F e_a,e_d)h(F e_b,e_c) \big)\\
				&\quad =2\big(\tr (E\circ F)\big)^2-2\tr\big((E\circ F)^2\big)
			\end{align*}
			Note that this part of the lemma does not use any properties of $E$ and $F$.
			\item The proof proceeds analogously, but there are several extra terms:
			\begin{align*}
				&\tr\big((K\obar_h K)\circ (L\obar_h L)\big)\\
				&\quad =\sum_{a,b,c,d}\epsilon_a\epsilon_b\epsilon_c\epsilon_d \\
				&\hspace{1.5cm}\big(h(K e_c,e_a)h(K e_d,e_b)-h(K e_c,e_b)h(K e_d,e_a)+2h(K e_c,e_d)h(Ke_a,e_b)\big)\\
				&\hspace{1.5cm}\big(h(L e_a,e_c)h(L e_b,e_d)-h(L e_a,e_d)h(L e_b,e_c)+2h(L e_a,e_b)h(Le_c,e_d)\big)\\
				&\quad =2\Big(\tr(K\circ L)\big)^2-2\tr\big((K\circ L)^2\big)
				+8\tr\big((L\circ K)^2\big)+4\big(\tr(K\circ L)\big)^2\Big)\\
				&\quad =6\big(\tr(K\circ L)\big)^2+6\tr\big((K\circ L)^2\big)
			\end{align*}
			Note that we used the fact that $K$ and $L$ are skew in arriving at the penultimate expression.
			\item Again, this is a straightforward computation:
			\begin{align*}
				&\tr\big((E\owedge_h F)\circ (K\obar_h K)\big)\\
				&\quad =\sum_{a,b,c,d}\epsilon_a\epsilon_b\epsilon_c\epsilon_d
				\big(h(E e_c,e_a)h(E e_d,e_b)-h(E e_c,e_b)h(E e_d,e_a)\big)\\
				&\hspace{1.5cm}\big(h(K e_a,e_c)h(K e_b,e_d)-h(K e_a,e_d)h(K e_b,e_c)+2h(K e_a,e_b)h(K e_c,e_d)\big)\\
				&\quad =2\big(\tr(E\circ K)\big)^2-2\tr\big((E\circ K)^2\big)-4\tr\big((E\circ K)^2\big)\\
				&\quad =2\big(\tr(E\circ K)\big)^2-6\tr\big((E\circ K)^2\big)
			\end{align*}
			where we used both self-adjointness of $E$ and skew-adjointness of $K$.
		\end{numberedlist}	
	\end{myproof}
	
	Now we must compute the relevant traces:
	
	\begin{lem}
		On $N_{n-1}$, the following trace identities hold for any non-negative integer $m$ and any $k\in \{1,2,3\}$:
		\begin{numberedlist}
			\item $\tr(\mc K^m)=4\Big((n-1)f_Z^m+ \frac{f_Z^{2m}}{f_\h^m}\Big)$.
			\item $\tr(\mc K^m I_k)=0$.
			\item $\tr(\mc K^m I_\h)=0$.
			\item $\tr(\mc K^m I_\h I_k)=0$.
		\end{numberedlist}
	\end{lem}
	\begin{myproof}
		Firstly, since $\mc K$ restricts to multiples of the identity on $\H Z$ and its orthogonal complement, it is easy to compute the traces of its powers:
		\begin{equation*}
			\tr(\mc K^m)=4\bigg((n-1)f_Z^m+\frac{f_Z^{2m}}{f_\h^m}\bigg)
		\end{equation*}
		In order to see that both $\tr(\mc K^m I_k)$ vanishes, we recall that $\mc K$ is self-adjoint with respect to $g$, while $I_k$ is skew, and that $\mc K$ and $I_k$ commute. This means that $\mc K^m I_k$ is skew, and has vanishing trace. The same argument applies to $\tr(\mc K^m I_\h)$. Finally, we consider $\tr(\mc K^m I_\h I_j)$. We know that $I_j=I_k I_l$, where $(j,k,l)$ is a cyclic permutation of $(1,2,3)$. This, combined with the fact that the endomorphisms all commute, yields 
		\begin{equation*}
			\tr(\mc K^m I_\h I_j)=\tr(\mc K^m I_\h I_k I_l)=\tr(I_l \mc K^m I_\h I_k)=\tr(\mc K^m I_\h I_l I_k)=-\tr(\mc K^m I_\h I_j)
		\end{equation*}
		whence $\tr(\mc K^m I_\h I_j)=0$.
	\end{myproof}
	
	Since, in this class of examples, $\mc K$, $I_\h$ and $I_k$ all commute and $I_\h^2=I_k^2=-\id$, any trace of a composition of these endomorphisms can be reduced to one of the four traces computed in the above lemma. That means that we are now in a position to compute the curvature norm.
	\begin{align*}
		32\norm{\mc R^\q}^2&=\tr(\id)^2-\tr(\id)+2\sum_k \big(\tr(I_k)^2-3\tr(I_k^2)\big)\\
		&\quad -\frac{2}{f_Zf_\h}\bigg(\tr(\mc K I_\h)^2-3\tr\big((\mc K I_\h)^2\big)
		+\sum_k\Big(\tr(\mc K I_\h I_k)^2-\tr\big((\mc K I_\h I_k)^2\big)\Big)\bigg)\\
		&\quad +3\sum_{j,k}\Big(\tr(I_j I_k)^2+\tr\big((I_jI_k)^2\big)\Big)\\
		&\quad -\frac{2}{f_Z f_\h}\bigg(3\sum_k \Big(\tr(I_k \mc K I_\h)^2+\tr\big((I_k \mc K I_\h)^2\big)\Big)\\
		&\hspace{1.9cm}
		+\sum_{j,k}\Big(\tr(I_k \mc K I_\h I_j)^2-3\tr\big((I_k \mc K I_\h I_j)^2\big)\Big)\bigg)\\
		&\quad +\frac{1}{f_Z^2 f_\h^2}\bigg(3\Big(\tr\big((\mc K I_\h)^2\big)^2+\tr\big((\mc K I_\h)^4\big)\Big)\\
		&\hspace{1.9cm}
		+2\sum_k\Big(\tr(\mc K I_\h \mc K I_\h I_k)^2-3\tr\big((\mc K I_\h \mc K I_\h I_k)^2\big)\Big)\bigg)\\
		&\quad +\frac{1}{f_Z^2 f_\h^2}\sum_{j,k}\Big(\tr(\mc K I_\h I_j \mc K I_\h I_k)^2
		-\tr\big((\mc K I_\h I_j \mc K I_\h I_k)^2\big)\Big)\\
		&=\tr(\id)^2-\tr(\id)+18\tr(\id)
		+3\sum_{j,k}\Big(\delta_{jk}\tr(\id)^2-(-1)^{\delta_{jk}}\tr(\id)\Big)\\
		&\quad -\frac{6}{f_Z f_\h}\Big(\tr(\mc K^2)
		-\sum_{j,k}(-1)^{\delta_{jk}}\tr(\mc K^2)\Big)\\
		&\quad +\frac{1}{f_Z^2 f_\h^2}\bigg(3\Big(\tr(\mc K^2)^2+\tr(\mc K^4)\Big)
		+18\tr(\mc K^4)\\
		&\hspace{1.9cm}+\sum_{j,k}\Big(\delta_{jk}\tr(\mc K^2)^2
		+(-1)^{\delta_{jk}}\tr(\mc K^4)\Big)\bigg)\\
		&=160n^2+32n+\frac{6}{f_Z^2f_\h^2}(\tr (\mc K^2)^2+4\tr(\mc K^4))\\
		&=32\Bigg(n(5n+1)
		+3\bigg((n-1)\frac{f_Z}{f_\h}+\frac{f_Z^3}{f_\h^3}\bigg)^2
		+3\bigg((n-1)\frac{f_Z^2}{f_\h^2}+\frac{f_Z^6}{f_\h^6}\bigg)\Bigg)
	\end{align*}
	This was precisely the claim.
\end{myproof}

Through the HK/QK correspondence, the function $2f_Z$ corresponds to a global coordinate function on $\bar N_{n-1}$, conventionally denoted by $\rho$ (and $2f_\h=-(\rho+2c)$). We have therefore just proven that $\norm{\mc R^\q}^2$ is a function only of $\rho$. This was to be expected in light of our results in~\cite{CST2020}, where we showed that the isometry group of $(\bar N_{n-1},g^c_\text{FS})$ acts transitively on the level sets of $\rho$. 

It follows from a short computation that $\norm{\mc R^\q}^2$ is an injective function of $\rho$ for any $c>0$, proving that the isometry group must preserve its level sets (for more details, see section 4 of~\cite{CST2020}). All this is summarized in the following theorem, which is also stated in the above-mentioned paper.

\begin{thm}
	The metrics $g_\text{FS}^{c>0}$ on $\bar N_n$ are of cohomogeneity one. In particular, for every non-negative integer $n$, the quaternionic K\"ahler manifold $(\bar N_n,g^{c>0}_\text{FS})$ is not a locally homogeneous space.\proofclear
\end{thm}

\end{document}